\documentclass[final,epsfig,amsfont]{article}
\usepackage{latexsym}
\usepackage[dvips]{graphicx}
\newcommand{\eeq}{\end{equation}}
\newcommand{\beq}{\begin{equation}}
\newcommand{\nuq}[1]{\label{#1} \eeq}

\newtheorem{conjecture}{Conjecture}

\newtheorem{remark}{Remark}
\newtheorem{lemma}{Lemma}
\newtheorem{corollary}{Corollary}

\newtheorem{theorem}{Theorem}
\newtheorem{proposition}{Proposition}

\begin{document}
\title{
Minkowski's question mark measure}

\author{
Giorgio Mantica \\
Center for Non-linear and Complex Systems, \\ Dipartimento di Scienza ed Alta Tecnologia, \\Universit\`a dell'~Insubria, 22100 Como, Italy \\ and \\ CNISM unit\`a di Como,
INFN sezione di Milano, \\
Istituto Nazionale di Alta Matematica, \\
Gruppo Nazionale per la Fisica Matematica.}

\date{}
\maketitle
\begin{abstract}
Minkowski's question mark function is the distribution function of a singular continuous measure: we study this measure from the point of view of logarithmic potential theory and orthogonal polynomials. We conjecture that it is regular, in the sense of Ullman--Stahl--Totik and moreover it belongs to a Nevai class: we provide numerical evidence of the validity of these conjectures. In addition, we study the zeros of its orthogonal polynomials and the associated Christoffel functions, for which asymptotic formulae are derived. As a by--product, we derive upper and lower bounds to the Hausdorff dimension of Minkowski's measure. Rigorous results and numerical techniques are based upon Iterated Function Systems composed of M\"obius maps.
\\
\\
{\em Keywords: Minkowski's question mark function; Orthogonal polynomials, Jacobi matrices; Regular Measures; M\"obius Iterated Function Systems; Nevai class; Gaussian integration; Christoffel functions.}\\
{\em MATH Subj. Class. 42C05; 47B36; 11A55; 11B57; 37D40 }
\end{abstract}

\section{Introduction: Minkowski's Q function and its singular measure}
\label{sec1}

\subsection{Definitions and goals of the paper}
Minkowski's question--mark function $Q(x)$ can be concisely defined---though not in the the most transparent way---by writing the point $x \in [0,1]$ in its continued fraction representation, $x = [n_1,n_2,\ldots,]$, by setting
$N_j(x) = \sum_{l=1}^j n_l$, and by defining $Q(x)$ as the sum of the series \cite{denjoy,salem}
\beq
 Q(x) = \sum_{j=1}^\infty (-1)^{j+1} 2^{-N_j(x)+1}.
\nuq{eq-minko1}
This function was originally constructed to map the rationals to the solutions of quadratic equations with rational coefficients in a continuous, order preserving way \cite{minko}, but it successively appeared that it has much wider implications in many fields of mathematics. 
A graph of $Q(x)$ is part of Figures \ref{fig-tetzet2} and \ref{fig-wei1b} below. It is remarkable 
that this graph can be seen as the attractor of an {\em Iterated Function System} (IFS) composed of M\"obius maps \cite{prlmob}, so that $Q(x)$ also belongs to the family of fractal interpolation functions \cite{ba3,ba4}.

In this paper we are interested in the singular-continuous measure $\mu$ for which $Q(x)$ is the distribution function:
\beq
 Q(x) = \int_0^x  d \mu.
\nuq{eq-minko2}
For short, in this paper we will call this measure the {\em Minkowski's measure} and we will reserve the symbol $\mu$ to denote it. It has been investigated in relation to singularity \cite{denjoy}, H\"older continuity \cite{salem}, the so--called thermodynamical formalism \cite{guzzi,kesse,kesse2,kinney}, the asymptotic behavior of its moments \cite{alau1,alau2} and of its Fourier transform  \cite{jordan,persson,salem,yakubo1,yakubo2}, so that many of its properties have been fully clarified.

To the contrary, the characterization of this measure from the point of view of orthogonal polynomials and logarithmic potential theory \cite{saff,stahl} is still
an open problem. This problem is the object of a recent publication by Dresse and Van Assche \cite{dresse}; by using a more powerful technique, which describes Minkowski's measure as the {\em balanced measure of an IFS} (to be defined below), we are in the position to correct their findings and to provide theoretical arguments and compelling numerical evidence in favor of the precise characterization of $\mu$ in terms of two principal conjectures:
\begin{conjecture}
\label{conj-reg}
Minkowski's  measure $\mu$ is regular in the sense of Ullman-Stahl-Totik.
\end{conjecture}
\begin{conjecture}
Minkowski's  measure $\mu$ belongs to the Nevai class $N(\frac{1}{4},\frac{1}{2}$).
\label{conj-nev}
\end{conjecture}

Let us briefly define the terms of these conjectures. They are linked to different asymptotic behaviors of the orthogonal polynomials $\{p_j(\mu;x)\}_{j \in {\bf N}}$ of Minkowski's measure $\mu$. Quite generally, given a positive Borel measure $\mu$ supported on a compact subset $E$ of the real axis, orthogonal polynomials are defined
by the relation $\int p_j(\mu;x) p_m(\mu;x) d \mu(x) = \delta_{jm}$, where $\delta_{jm}$ is the Kronecker delta. They satisfy the three-terms recurrence relation
\begin{equation}
\label{nor2}
   x p_j (\mu;x) = a_{j+1} p_{j+1}(\mu;x)
   + b_j p_j(\mu;x)  + a_{j} p_{j-1}(\mu;x),
\end{equation}
initialized by $a_0=0$ and
$p_{-1}(\mu;s) = 0$, $p_0(\mu;s) = 1$, that can be encoded
in the Jacobi matrix $J(\mu)$:
\beq
   J(\mu) :=
           \left(   \begin{array}{ccccc}  b_0  & a_1 &      &      &       \cr
                               a_1  & b_1 & a_2  &      &       \cr
                      &   \ddots    &    \ddots      & \ddots &\cr
                               \end{array} \right) .
 \nuq{jame}
For compact support $E$ the moment problem is determined \cite{ach} and the matrix $J(\mu)$
is in one--to--one relation with the measure $\mu$.
At times, it is also important to consider {\em monic} orthogonal polynomials, $P_j(\mu;x)$, normalized so that $P_j(\mu;x) = x^j + q_{j-1}(x)$, in which $q_{j-1}$ is a polynomial of degree $j-1$. Of course, the two differ by a constant factor: $p_j = \gamma_j P_j$ and from eq. (\ref{nor2}) it easily follows that $1/\gamma_j = \prod_{l=1}^j a_l$.

A well established theory, started by the classical works of Geronimus and Szeg\"o \cite{gero,szego},
classifies different asymptotic behaviors of these polynomials, which can be summarized as follows, in the case of a measure whose support is $[0,1]$:

{\bf Root asymptotics:} For $z$ not in $[0,1]$, when the degree $j$ tends to infinity, $p_j(\mu;z)^{1/j}$ tends to $\phi(z)=2z-1+2\sqrt{z^2-z}$, the function that maps $[0,1]$ to the unit circle in the complex plane; equivalently, $\gamma_j^{-1/j} \rightarrow \mbox{Cap} ([0,1]) = \frac{1}{4}$.

{\bf Ratio asymptotics:} For $z$ not in $[0,1]$, when the degree $j$ tends to infinity, $p_{j+1}(\mu;z)/p_j(\mu;z)$ tends to $\phi(z)$; equivalently, $a_j=\gamma_{j-1}/\gamma_j \rightarrow \frac{1}{4}$, $b_j \rightarrow 0.$

{\bf Power asymptotics:} For $z$ not in $[0,1]$, when the degree $j$ tends to infinity, $p_{j}(\mu;z)/\phi(z)^j$ tends to a function $g(z)$; equivalently, $\gamma_{j}/4^j \rightarrow \gamma>0$.

Convergence in the above is understood to be uniform for $z$ in compact sets in the complement of $[0,1]$ in the complex plane.
Observe from the start that Minkowski's singular continuous measure $\mu$ can not fulfill power asymptotics, for otherwise it would belong to Szeg\"o class, {\em i.e.} it would be absolutely continuous, with a density satisfying the well known Szeg\"o condition. Therefore, if validated, Conjectures \ref{conj-reg} (which is equivalent to root asymptotics) and \ref{conj-nev} (equivalent to ratio asymptotics) would grant the strongest possible regularity to which $\mu$ can aspire, in the above hierarchy.

It is remarkable that these asymptotics have a simple, equivalent definition in terms of the Jacobi matrix of $\mu$. Therefore, the first goal of this paper is its precise computation---for which we will employ an algorithm that we have introduced in \cite{mobius}---and the numerical verification of the above conjectures.

The second main goal of the paper, linked to the first, is the investigation of three sequences of discrete measures that are related to $\mu$ and of their limits.
In this paper, convergence in the space ${\cal M}([0,1])$ of regular Borel probability measures on $[0,1]$ will always be understood in weak star sense, {\em weakly} for short, even when not specified: a sequence $\mu_j$ converges to $\mu$ if and only if the integrals of any continuous function $f$ with respect to $\mu_j$ tend to the integral of $f$ with respect to $\mu$.
The measures under investigation are supported on the zeros of the orthogonal polynomials of $\mu$. We will use the notation $\zeta^j_l$ to denote these zeros:
\beq
p_j(\mu;\zeta^j_l) = 0, \; l=1,\ldots,j.
\nuq{eq-zero0}
These measures are: the so--called measure of the zeros, the measure underlying Gaussian integration, the measure linked to root asymptotics.  Let us briefly review their definition and main properties.

The {\em measure of zeros} $\nu_j$ is obtained by placing atomic measures of equal weight at the zeros of the orthogonal polynomials:
\beq
   \nu_j = \frac{1}{j} \sum_{l=1}^j  \delta_{\zeta^j_l}.
 \nuq{eq-nu3}
Here and in the following, $\delta_x$ is the unit mass atomic measure supported at $x$.
The measure $\nu_j$ is termed the density of states in the physical literature. It can be used to define a further regularity property.
Suppose that $E$, the support of $\mu$, has positive capacity and that $\nu_j$ converges weakly to the {\em equilibrium measure $\nu_E$} on this support \cite{saff} when the polynomial degree $j$ tends to infinity. Then, the orthogonal polynomials of $\mu$ are said to have {\em regular asymptotic zero distribution}: this is definition 3.1.3 of \cite{stahl}, that is a reference also for the following result.
Observe that the support of Minkowski's measure is the interval $E=[0,1]$, which has empty interior (as a set in the complex plane) and positive capacity. In such a case, {\em regularity of the zero distribution is equivalent to regularity of the measure $\mu$} (Theorem 3.1.4 \cite{stahl}), {\em which is also equivalent to root asymptotics} (Theorems 3.1.1 and 3.2.1 \cite{stahl}).

The {\em measure underlying Gaussian integration} $\mu_j$ is obtained replacing the equal weights in eq. (\ref{eq-nu3}) with the so--called {\em Christoffel numbers} $w^j_l$ (that we shall also call {\em Gaussian weights} and are defined in  eq. (\ref{eq-nu2}) below):
\beq
   \mu_j =  \sum_{l=1}^j w^j_l \delta_{\zeta^j_l}.
 \nuq{eq-nu1}
At difference with $\nu_j$, the sequence $\mu_j$ always converge: it is well known that integrals of polynomials up to degree $2j-1$ with respect to $\mu$ and $\mu_j$ coincide, so that $\mu_j$ tends to $\mu$ weakly.

The {\em measure linked to root asymptotics} $\sigma_j$ is the third type of measure that we will study in the paper. It is defined via yet another choice of the weights placed at the location of the zeros:
\begin{equation}
\sigma_j = \sum_{l=1}^j w^j_l  p_{j-1}^2(\mu;\zeta^j_l) \delta_{\zeta^j_l}.
\label{eq-sig02b}
\end{equation}
Weak star convergence of the sequence $\sigma_j$ is equivalent to the fact that $\mu$ belongs to the {\em Nevai class of measures} $N(a,b)$ \cite{walter2}. Measures in this class are those for which diagonal and out-diagonal Jacobi matrix entries tend to a limit as $j$ tends to infinity:
\beq
   \lim_{j \rightarrow \infty} (a_j,b_j) = (a,b).
 \nuq{eq-nu5}
In our case, because of symmetry, $b_j=1/2$ for all values of $j$.
In turn, Nevai class is related, as seen above, to ratio asymptotics.

Clearly, since the Nevai class is a subset of that of regular measures, proving that Minkowski's measure $\mu$ belongs to this class also implies that it is regular, so that we could have focused only on this last problem. Nonetheless, we will describe both conjectures separately and in order of difficulty, because much is to be learned in each step.

\subsection{Organization of the paper and summary of results}

The peculiarities of the Minskowski measure are best revealed when seeing it as a balanced measure of a M\"obius Iterated Function System.
In the next section we review {\em two} IFS, composed of M\"obius maps, that can be associated to Minkowski's question mark function  \cite{prlmob,mobius} and we briefly comment on their relations with the theory of dynamical systems.
We show how they can be used to construct different sequences of probability measures that converge weakly to Minkowski's measure $\mu$.

In Section \ref{sec-mobalg} the previous theory is translated into two algorithms for the computation of the Jacobi matrix $J(\mu)$. The first is based upon the Jacobi matrix of a finite sum of atomic measures, which approximates $\mu$. We show that it is largely inefficient, when applied to the present case. This also explains the lack of convergence observed in Section 2 of \cite{dresse}. The second is the technique proposed in \cite{mobius}: in two subsections we test it for speed of convergence and precision. These tests are mandatory, since the numerical results of this paper rely on the computed entries of the Jacobi matrix $J(\mu)$. We justify experimentally the claim that this algorithm, when run in double precision Fortran on a standard desktop computer, can provide Jacobi matrices of rank about 60,000 with a controlled error.

The main part of the paper presents the investigation of the conjectures mentioned before, as well as further conjectures and finer details of Minkowski's measure and of its orthogonal polynomials. We start in Section \ref{sec-mink} by analyzing the regularity of the measure $\mu$, {\em i.e.} Conjecture \ref{conj-reg}, via convergence of the sequence $\gamma_j^{-1/j}$ to the value $\frac{1}{4}$. We observe that this convergence holds and it is of power--law type: there exist two positive constants $A$ and $B$ so that $|\gamma_j^{-1/j}-\frac{1}{4}|\sim A j^{-B}$. The capital letters $A$ and $B$ will indicate throughout the paper different pairs of constants that appear in power-law behaviors.

Next, we study the {\em regular asymptotic distribution of the zeros}, as defined above. In Section \ref{sec-zeroreg}  we provide numerical evidence of a stronger property, expressed by \\ {\bf Conjecture \ref{conjzer}}: {\em The zeros of Minkowski's polynomials converge uniformly to the zeros of the Chebyshev polynomials of the same order, as $j$ tends to infinity}. \\
In Section \ref{sec-discre} regularity of the distribution of the zeros is established via the analysis of the so--called {\em discrepancy} \cite{andriev3,andriev4,andriev5} between the equilibrium measure $\nu_E$ and the sequence of measures $\nu_j$.  Our findings are summarized in \\ {\bf Conjecture \ref{conjdiscre}}: {\em As $j$ tends to infinity, the discrepancy $D(\nu_E,\nu_j)$ tends to zero and convergence is of power--law type}.

In the central Section \ref{sec-weights} we examine the measure of Gaussian integration $\mu_j$ in eq. (\ref{eq-nu1}). This detailed analysis best reveals the characteristic features of Minkowski's measure $\mu$. Together with Gaussian nodes ({\em i.e.} the zeros of the orthogonal polynomials) that have been studied in the previous sections, we here investigate Christoffel numbers. They can be numerically computed in different ways \cite{gautschi,gautchr}, the most common being the Golub--Welsh algorithm \cite{gene}. To the contrary, this latter fails in the case of Minkowski's measure, but in its stead we profitably use a formula due to Shohat \cite{shoha} (see eq. (\ref{eq-nu2}) below) based on the reproducing kernel of the orthogonal polynomials.

The reason behind the failure of Golub--Welsh in this case is a distinctive feature of Minkowski's measure $\mu$: the presence of extremely small Christoffel numbers even at moderate polynomial orders $j$---that become of the order of $10^{-1000}$ for $j\simeq 60,000$.
We describe a reliable technique for the computation of the atoms composing the measure $\mu_j$, which can be appropriately dubbed {\em mathematical neutrinos}, for their elusive mass. The technique is inspired by a common practice in the calculation of Lyapunov exponents. A practical result of this computation is the derivation of rapidly converging upper and lower bounds to the value of the Hausdorff dimension of the measure $\mu$, in Sect. \ref{sec-hausdim}.

The theoretical analysis of Gaussian weights--Christoffel functions is carried out to a large detail in Sections \ref{sec-slip} and \ref{sec-numconf}, Propositions \ref{prop-fundam} and \ref{cor-chr}, which can be briefly and informally summarized as follows:
\\
{\bf Proposition \ref{prop-fundam}}: {\em For sufficiently large $j$, the logarithm of the average amplitude of weights $w^j_{l}$, when $\zeta^j_{l}$ is in the neighborhood of the rational point $\frac{1}{q}$, is bounded between two explicit functions reported in eq. (\ref{eq-broc40}).} \\
{\bf Proposition \ref{cor-chr}}: {\em The logarithm of the Christoffel functions $\lambda_j (x)$ has an asymptotic singular behavior of the kind $\log(\lambda_j (x)) \sim - \log(2)/(q^2 |x-\frac{1}{q}|) - \log(j)$, for $x$ in the neighborhood of $\frac{1}{q}$.} \\
These results are a consequence of regularity of the measure $\mu$ and of the non--analytic behavior of $Q$ at rational values, described in {\bf Lemma \ref{lem-farey}}, Section \ref{sec-slip}.
The theory presented in Section \ref{sec-weights} contains much more information, which is discussed at length with the aid of detailed figures and fits with the hierarchical structure of the so-called {\em Farey tree}, defined in {\bf Remark \ref{rem-farey}}, Section \ref{sec-suppo}. We perform the analysis around the Farey points $\frac{1}{q}$, $q$ integer, but the techniques employed in Lemma \ref{lem-farey} permit to treat the case of any rational point.

Finally, in Section \ref{sec-nevai}, we try to assess whether Minkowski's question mark measure belongs to the Nevai class $N(\frac{1}{4},\frac{1}{2})$. We do this in two ways. The first is the direct verification of convergence of Jacobi matrix elements. Our data reveal that the out--diagonal entries tend to the limit value one quarter at a slow pace, which nonetheless seems to be of power--law type. This yields a stronger version of Conjecture \ref{conj-nev}: \\
{\bf Conjecture \ref{conj-nevpower}}: {\em Convergence of the Jacobi matrix elements of Minkowski's measure is of power--law type: there exist two positive constants $A,B$ such that $|a_j-\frac{1}{4}| < A j^{-B}$}.
\\
In the process, we also show that the computed matrix elements are consistent with the theoretical result that stronger regularity properties, such as Szeg\"o asymptotics, do not hold.
The second method is via the analysis of the measures $\sigma_j$ defined above in eq. (\ref{eq-sig02b}), which are associated with ratio asymptotics: by showing numerically that the Hutchinson distance between $\sigma_j$ and the expected limit measure $\sigma_E$ vanishes when $j$ tends to infinity we provide further evidence in favor of enlisting $\mu$ in the Nevai class: \\
{\bf Conjecture \ref{conj-nevsigma}}: {\em
Convergence of the sequence of measures $\sigma_j$ to $\sigma_E$ is of power--law type: there exist two positive constants $A,B$ such that the Hutchinson distance verifies
$ d(\sigma_j,\sigma_E) < A j^{-B}$.}

In the Conclusion we briefly comment on the relevance of these results and on possible extensions of this analysis.

\section{M\"obius Iterated Function Systems}
\label{sec-mobifs}

Minkowski's question mark function also appears in dynamical systems: for instance, it makes the conjugation of the two dynamical systems generated by the tent and Farey maps on the unit interval \cite{rugh,mirko,isola}. But a second aspect is more important for the present work.

In 1988 Gutzwiller and Mandelbrot \cite{guzzi} studied {\em coding functions} for the motion of billiards in hyperbolic space \cite{series}. Because of the negative curvature of these spaces, these trajectories are highly unstable and indeed they provide a principal model of chaotic system \cite{arnold}. Gutzwiller's motivation was to find a purely abstract analogue of the coding function of the anisotropic Kepler problem \cite{guzzi-ke}, another strongly chaotic system.
By its very nature, a coding function encodes in its value the infinite symbolic history of a trajectory. Without entering in the details of the definition of this function, it is enough to remark that, in the case of a particular triangular billiard, the Gutzwiller--Mandelbrot coding function can be re-written in the form (\ref{eq-minko1}) and hence it coincides with Minkowski's $Q$ function. By employing the symmetry properties of such hyperbolic billiards, we were able to show that this function---as well as the coding functions of more general hyperbolic billiards and also of the anisotropic Kepler system---can be described by {Iterated Function Systems} composed of M\"obius maps \cite{prlmob}.

Recall that a hyperbolic IFS \cite{hut,dem,ba2}
is a collection of contractive maps $\Phi_i : X \rightarrow
X$, $i = 1, \ldots, M$, on a compact metric space $X$. For these systems, there
exists a unique compact set ${\mathcal A}$, called the {\em attractor} of the IFS, that solves the equation
 \begin{equation}
 \label{attra}
    {\mathcal A}=\bigcup_{i=1,\ldots ,M}\;\Phi_i({\mathcal A}) :=  \Phi ({\cal A}).
 \end{equation}

The first IFS associated with Minkowski's function $Q$ acts on the unit square, $X=[0,1]^2$, and has the {\em graph} of $Q$ as attractor.
This IFS can be constructed as follows.
Let the maps $M_i$ and $P_i$, $i=1,2$ be defined as
\begin{equation}
\label{eq-mink1}
    \begin{array}{ll}
    M_1(x) = \frac{x}{1+x}, & P_1(y) = \frac{y}{2},\\
    M_2(x) = \frac{1}{2-x}, & P_2(y) = \frac{y+1}{2}
    \end{array}
\end{equation}
and let us write
\begin{equation}
\label{eq-mink11}
    \Phi_i(x,y) = (M_i(x),P_i(y)), \; i = 1,2.
\end{equation}
Then, the following results hold \cite{prlmob}:
\begin{theorem}
Let $\{\Phi_i(x,y), \; i = 1,2 \}$ be the IFS on the unit square defined in eqs. (\ref{eq-mink1}),(\ref{eq-mink11}). Its attractor is
${\mathcal A} = \{ (x,Q(x)), x \in [0,1] \}$ and ${\mathcal A} = \lim_{n \rightarrow \infty} \Phi^n (K)$ in the Hausdorff metric, where $K$ is any non--empty compact set in $[0,1]^2$.
\end{theorem}

\begin{corollary}
It is convenient to take $K = \{ (x,x), x \in [0,1] \}$, the graph of the identity function: then, for any $n$, $\Phi^n (K)$ is the graph of a function $Q_n(x)$ that tends to Minkowski's function $Q(x)$ uniformly when $n$ tends to infinity.  Thanks to the above, one also proves the symmetries
\begin{equation}
\label{eq-mink11b}
    Q(M_i(x)) = P_i(Q(x)), \; i = 1,2.
\end{equation}
\label{rem-one}
\end{corollary}
This procedure was used to produce the graph of $Q(x)$ in Fig. 2 of \cite{prlmob}, Fig. 1 of \cite{dresse} and Figures \ref{fig-tetzet2} and \ref{fig-wei1b} herein.

The second IFS we consider is associated with Minkowski's measure $\mu$ and it acts in $X=[0,1]$.
It is composed of the two M\"obius maps $\phi_i = M_i$ just introduced. Clearly, using now the lower case maps $\phi_i$ to define the IFS in eq. (\ref{attra}) implies that ${\mathcal A} = [0,1]$, but our interest at this point is not in the attractor.

Rather, recall that a {\em balanced measure} can be uniquely associated with a hyperbolic IFS $\{\phi_i : [0,1] \rightarrow
[0,1]$, $i = 1, \ldots, M\}$, through the choice of a set of probabilities $\rho_i>0$, $i = 1, \ldots, M$: $\sum_i \rho_i = 1$ \cite{ba2,dem,hut}. The theory has been extended by Mendivil to transformations that are contractive on average \cite{mendiv}. A transformation $T^*$ on the space ${\cal M}([0,1])$ of regular Borel probability measures on $[0,1]$ can be defined as follows: $T^* \eta $ is the unique measure that verifies the equation
\begin{equation}
\label{bala}
   \int f \; d (T^* \eta) \; =
   \sum_{i=1}^{M}
   \; \rho_{i}
   \; \int   \;
 (f \circ \phi_{i}) \; d\eta ,
\end{equation}
for any continuous function $f$. $T^*$ is also known as the Perron--Frobenius operator. This operator is contractive in the Hutchinson--Wesserstein--Kantorovich metrics (which is a metric that entails weak star convergence) under which ${\cal M}([0,1])$ is a complete metric space: if $\eta$ and $\theta$ belong to ${\cal M}([0,1])$, their distance can be defined as
\beq
d(\eta,\theta) = \sup \{\int f d\eta - \int f d\theta \},
\nuq{eq-hutdis}
where the supremum is taken over the set of $\mbox{Lip}_1$ functions.
The general theory applied to our case provides the following results \cite{prlmob,mobius}:

\begin{theorem}
\label{teo-two}
Let $\{M_i(x), \; i = 1,2 \}$ be the IFS maps defined in eq. (\ref{eq-mink1}), let $\rho_1=\rho_2=\frac{1}{2}$ and let $T^*$ be the Perron--Frobenius operator defined in eq. (\ref{bala}) with $\phi_i=M_i$. The balanced measure of this IFS coincides with Minkowski's measure $\mu$:
\begin{equation}
\label{bala2}
   T^* \mu = \mu.
\end{equation}
Moreover, for any $\eta \in {\cal M}([0,1])$, the sequence of measures $\eta_n$ converges weakly in ${\cal M}([0,1])$ to Minkowski's measure $\mu$
\begin{equation}
\label{eq-rec1}
   \eta_n := (T^*)^n\eta \rightarrow \mu \;\; \mbox{  when } n \rightarrow \infty.
\end{equation}
\end{theorem}

\begin{corollary}
On a par with Corollary \ref{rem-one}, one can take $\eta$ as the uniform measure in $[0,1]$ and compute Minkowski's measure $\mu$ as the limit of $(T^*)^n \eta$.
\end{corollary}

\section{The M\"obius IFS algorithm}
\label{sec-mobalg}

\subsection{Supportive theory and brief description}
\label{sec-suppo}
Let us now describe how the Jacobi matrix of $\mu$ can be computed. The full supportive theory has indeed been presented in \cite{mobius}, so that only a brief summary of it is reviewed here.
The first immediate observation is that, being the moment problem for $\eta$ and $\eta_n$ determinate, weak convergence of $\eta_n$ implies strong convergence of the corresponding Jacobi matrices (Theorem 1 in \cite{mobius}):
\begin{theorem}
\label{teo-3}
Let $J(\eta_n)$ be the Jacobi matrix of $\eta_n$, defined as in eq. (\ref{eq-rec1})  and let $(a_j(\eta_n),b_j(\eta_n))$ be its entries. Then,
\begin{equation}
\label{eq-num01}
   (a_j(\eta_n),b_j(\eta_n)) \rightarrow (a_j(\mu),b_j(\mu)) \;\; \mbox{  when } n \rightarrow \infty
\end{equation}
for any value of $j$.
\end{theorem}

This result can be used in two ways.
Firstly, one can compute the sequence of measures $\eta_n$ and successively $J(\eta_n)$. This is similar to the approach followed by Dresse and Van Assche in Sect. 2 of \cite{dresse}: they have chosen a weakly convergent sequence of discrete measures $\eta_n$, composed of $2^n+1$ atoms, and computed by the discretized Stieltjes algorithm the corresponding Jacobi matrices.
In the IFS setting, one can take an initial measure $\eta$ composed of a finite number $L$ of atoms, like {\em e.g.} $\eta=\delta_{\frac{1}{2}}$ or $\eta=\frac{1}{2} ( \delta_0+\delta_1)$. Therefore, if
\beq
  \eta_n = \sum_{l=1}^L c^n_l \delta_{x^n_l},
  \nuq{eq-num02}
the action of the Perron--Frobenius operator $T^*$ yields
\beq
 \eta_{n+1}= T^* \eta_n = \frac{1}{2} \sum_{i=1}^2 \sum_{l=1}^L c^n_l \delta_{M_i(x^n_l)}.
  \nuq{eq-num02b}
The new measure is still composed of a finite number of atoms, $2L$, at the positions $M_i(x^n_l)$ and with the new weights $c^n_l/2$.
\begin{remark}
If $\eta=\delta_{\frac{1}{2}}$ the atoms of $\eta_n$ are located at  $x^n_l$, which are the nodes of level $n$ of the Farey tree, rooted at $x^0_1=\frac{1}{2}$.  This tree can therefore be defined recursively via eq. (\ref{eq-num02b}).
\label{rem-farey}
\end{remark}
The Farey tree obviously plays a major role in the theory of the Minkowski measure, as will appear below.

To compute the Jacobi matrix of the discrete measure $\eta_n$, one can use the discretized Stieltjes algorithm \cite{dresse}, or better Gragg and Harrod's algorithm \cite{gragg}, which  is more stable and can be enhanced to treat large sets of atoms \cite{arx}. Nonetheless, this approach is only helpful if relatively small values of $n$ are sufficient to yield a significant number of Jacobi matrix entries (at convergence), because of the geometrical increase with $n$ of the computational complexity. As remarked in \cite{dresse} and proven below we are in the opposite situation, which jeopardizes this approach.

Therefore, a different approach is needed. While Dresse and Van Assche considered a Chebyshev method based on ordinary moments, we employ the technique that we developed in \cite{mobius}. We start from equations (\ref{bala}) and (\ref{eq-rec1}), but we read them differently, as a transformation of Jacobi matrices into Jacobi matrices, in which the measures $\eta_n$ {\em neither appear nor have to be computed}.
This magic is effected in two steps \cite{mobius}.

{\bf 1:} Using the spectral theorem for $J(\eta)$, one defines $M_i(J(\eta))$ for the two M\"obius maps in eq. (\ref{eq-mink1}): this is the content of {Theorem 4} in \cite{mobius}, that also shows that numerical computation of $M_i(J(\eta))$ requires the solution of a tridiagonal linear system.

{\bf 2:} The two Jacobi matrices $M_i(J(\eta))$, for $i=1,2$ are in one--to--one correspondence with two measures, so that we use the algorithm of Elhay--Golub--Kautsky \cite{elhay} to compute the Jacobi matrix of their arithmetical average, which is precisely $J(T^*(\eta))$.

From an abstract point of view, the combination of {\bf 1} and {\bf 2} defines a map in the set of Jacobi matrices of measures in ${\cal M}([0,1])$: we call it
${\cal T} : {\cal J}([0,1]) \rightarrow  {\cal J}([0,1])$,
\beq
  {\cal T} ( J(\eta)) = J(T^*(\eta)),
  \nuq{eq-num03}
so that one can compute iteratively the left hand side $(a_j(\eta_n),b_j(\eta_n))$ in Theorem \ref{teo-3}, eq. (\ref{eq-num01}).

Original numerical experiments \cite{mobius} showed convergence of Jacobi matrix elements up to $j=4,000$ (see Fig. 3 in Sect. 6 of \cite{mobius}, where following Gutzwiller--Mandelbrot \cite{guzzi} Minkowski's $Q$ function is called the {\em slippery devil's staircase}). In this paper we need to push the technique of \cite{mobius} to much larger orders, to resolve the conjectures presented in the Introduction.
To do this without peril, we must first assess convergence and reliability of the M\"obius IFS algorithm.


\subsection{Speed of convergence}
\begin{figure}
\centerline{\includegraphics[width=.6\textwidth, angle = -90]{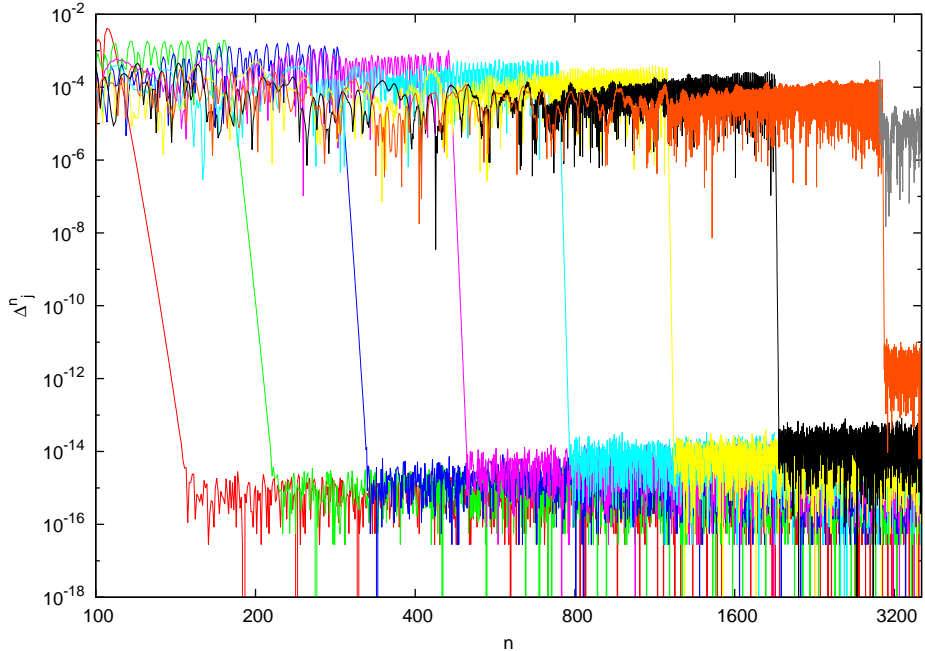}}
\caption{Difference $\Delta^n_j = |a^n_j-a^{n-1}_j|$ versus generation numbers $n$, for $j=125, 250, 500, 1000, 2000, 4000, 8000, 16000, 32000$. Different curves can be ordered according to the location of their sudden drop, increasing with $j$ from left to right.}
\label{erra2}
\end{figure}

First of all, we test the speed of convergence in eq. (\ref{eq-num01}). We start from the Jacobi matrix of the uniform measure on $[0,1]$ (that is, from the recurrence relation of rescaled Legendre polynomials). Of course, implementation of (\ref{eq-num03}) requires a truncation of the Jacobi matrices involved: we fix a large truncation, $j_{\max}=40,000$, and we compute $\Delta^n_j = |a^n_j-a^{n-1}_j|$, that is, the variation of an element of the Jacobi matrix from the $n$-th iteration of the algorithm to the next. These differences are plotted versus $n$ in Figure \ref{erra2}, for geometrically increasing values of $j$. We observe that the oscillations $\Delta^n_j$ suddenly drop by orders of magnitude. This is a clear indication that convergence has been reached in a finite subspace.

\begin{figure}
\centerline{\includegraphics[width=.6\textwidth, angle = -90]{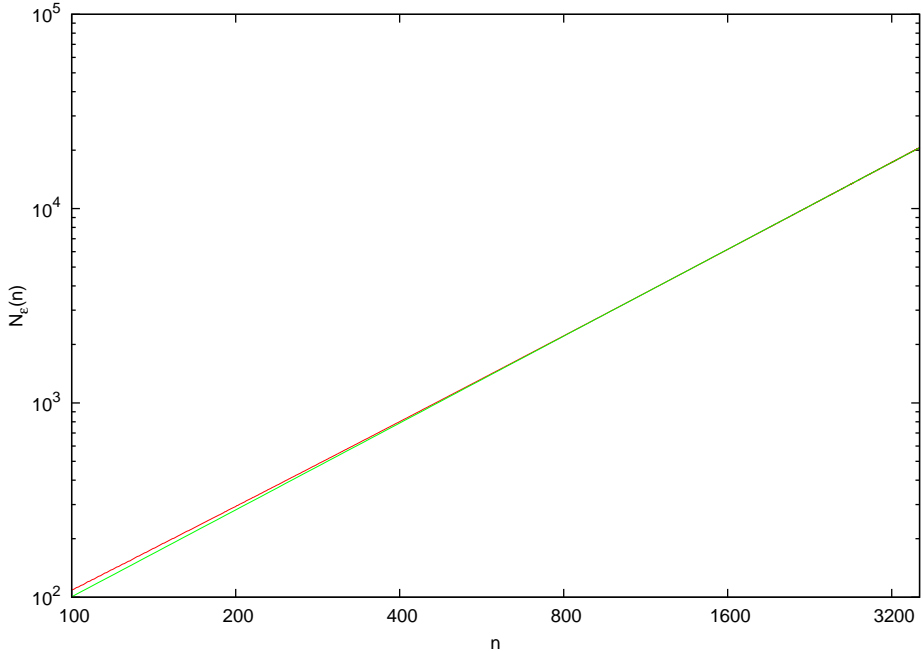}}
\caption{Range of numerical convergence $N_\epsilon(n)$ (red) versus $n$ in the case of Fig. \ref{erra2} with $\epsilon = 10^{-3}$. The fitting line (green) is a power--law, $N_\epsilon(n) = A n^{B}$ with $A=.11$ and $B=1.48$.}
\label{range1}
\end{figure}

To verify this assumption and to quantify the number of iterations required, we define the {\em rank of the Jacobi matrix at numerical convergence within a threshold $\epsilon$ after $n$ iterations} as
\beq
  N_\epsilon(n) := \max \{N \mbox{  s.t. } \sum_{j=1}^N |a^n_j-a^{n+1}_j| \leq \epsilon \}.
  \label{eq-neps3}
  \end{equation}
In fact, according to eq. (\ref{eq-neps3}) the vector of matrix elements $a_j$ with $1 \leq j \leq N_\epsilon(n)$ vary less than $\epsilon$ in $L_1$ norm when one further iteration of the algorithm is effected after the $n$-th.
This quantity is plotted versus $n$ in Figure \ref{range1}: one observes a regular  increase of $N_\epsilon(n)$ with $n$ which indicates three facts. The first is that convergence has been stably reached in the range $1 \leq j \leq N_\epsilon(n)$ independently of truncation. The second is that the number of converged components increases at a power--law rate: $N_\epsilon(n) \sim A n^B$. The third is that, while this increase is fully manageable in the M\"obius algorithm, it prevents the usage of discrete measures, as mentioned in the previous subsection: in fact,
assuming the same rate of convergence and required precision, a number of atoms of the order of the exponential of $N^{1/B}$ would be needed in eqs. (\ref{eq-num02}),(\ref{eq-num02b}) to compute $N$ Jacobi matrix elements.

The $L_1$ norm employed in eq. (\ref{eq-neps3}) entails a rather cogent test since it also implies a bound on the infinity norm. We can now compare our results with the table published in \cite{dresse}. These latter data were obtained via the Chebyshev algorithm from ordinary moments, programmed in Maple with 400 digits arithmetics, to cope with the numerical instability of the technique. It can be assumed that all digits reported are correct.  We have therefore taken the square roots of the tabulated $a_j^2$ coefficients and retained the same number of significant digits, 20. We compare these values with the values provided by M\"obius IFS algorithm, coded in double precision Fortran on a common desktop computer, reported with 15 significant digits. We observe almost perfect agreement, the difference between the two data sets being less than $10^{-15}$ in all cases.
This fact is consistent with the data of Figure \ref{erra2}: the values of $\Delta^n_j$, for $j=125$ (first curve, red) drop to a value of about $10^{-15}$ after convergence, which takes place at about $n=150$.

\begin{table}
\label{tabba}
\centering
\begin{tabular}{|r|l|l|}
  \hline
  $j$ & Algorithm &  $a_j$  \\
  \hline
  $1$ & M\"obius & 0.202302932329981 \\
  $1$ & Chebyshev &  0.20230293232998066551  \\
 $10$ & M\"obius& 0.215070562228327 \\
 $10$ & Chebyshev &  0.21507056222832743181  \\
 $20$ & M\"obius& 0.224458577806858 \\
 $20$ & Chebyshev &  0.22445857780685732313  \\
 $30$ & M\"obius& 0.221516521450380 \\
 $30$ & Chebyshev &  0.22151652145038000730  \\
 $40$ &  M\"obius& 0.236423204888560 \\
 $40$ & Chebyshev &  0.23642320488855968894 \\
  \hline
\end{tabular}
\caption{Jacobi matrix elements $a_j$, for $j=1,10,20,30,40$ from the M\"obius algorithm and from Table 2 in \cite{dresse} (Chebyshev algorithm).}
\end{table}

Unfortunately, the comparison in Table \ref{tabba} cannot be exploited much further, because the Chebyshev technique of \cite{dresse} cannot be easily extended beyond $j=40$. We have therefore to develop different tests that rely uniquely on our data.
A first indication follows again from Figure \ref{erra2}, in which we observe that, after convergence, the values of $\Delta^n_j$ oscillate around an average value that {\em increases} with $j$. This increase is a clear sign of loss of precision that we now investigate.

\subsection{Precision of the truncated, fixed point Jacobi matrix}

\begin{figure}
\centerline{\includegraphics[width=.6\textwidth, angle = -90]{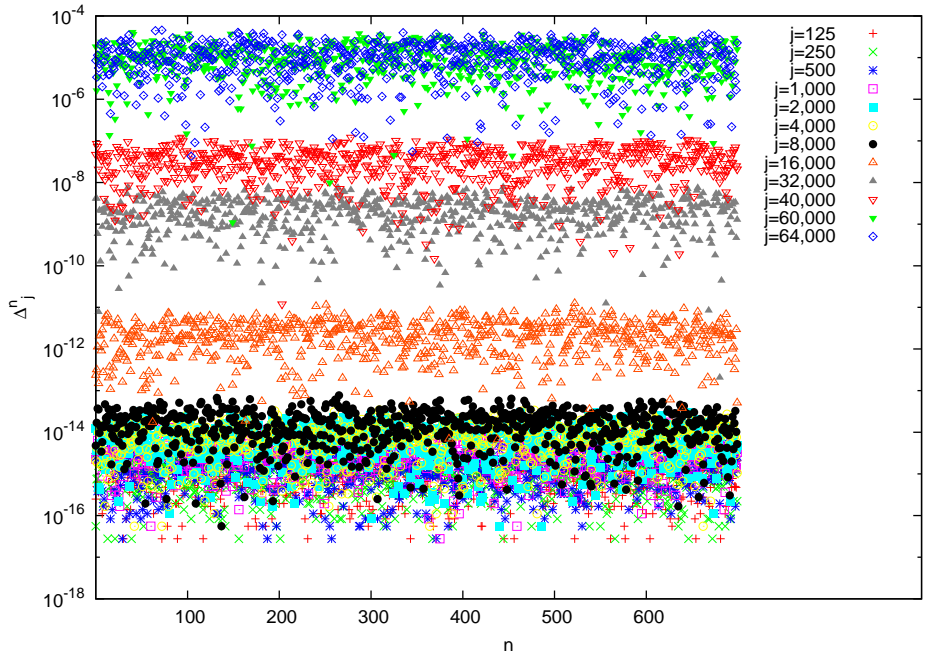}}
\caption{Differences $\Delta^n_j = |a^n_j-a^{n-1}_j|$ versus iteration $n$, for the values of $j$ reported. Initial Jacobi matrix is $J_{\mbox{\tiny num}}$, the limit matrix obtained after 8,000 iterations of the algorithm. }
\label{oscilla}
\end{figure}

The M\"obius IFS algorithm is an iterative technique that converges to a fixed point. The details of such convergence clearly depend on the initial point of the iteration, {\em i.e.} the Jacobi matrix input to the procedure. In our case, we have used the Jacobi matrix of Legendre polynomials on $[0,1]$. In any case, as observed in Figures \ref{erra2} and \ref{range1}, after a certain number of iterations, numerical convergence is reached in a finite subspace. To investigate this phenomenon we set the size of Jacobi matrices to the largest experimented value, $j_{\max}=64,000$ and we run the algorithm for $N=8,000$ iterations.
We observe that {all} matrix elements in the range $1 \leq j \leq 64,000$ reach numerical convergence:
an approximation, $J_{\mbox{\tiny num}}(\mu)$, of the fixed point Jacobi matrix,
\beq
\begin{array}{lcl}
  {\cal T} ( J(\mu) ) =  J(\mu),
  \end{array}
  \nuq{eq-num04}
has been evaluated. To verify the precision of this approximation, we now take as {\em input} the computed Jacobi matrix $J_{\mbox{\tiny num}}(\mu)$ and we further run the algorithm, that is, we act on $J_{\mbox{\tiny num}}(\mu)$ by  ${\cal T}$ for another $N=700$ iterations. Of course, the values of $a_j(n)$ are not constant, but fluctuate, due to numerical errors and finite truncation. In Figure \ref{oscilla} we plot again $\Delta^n_j = |a^n_j-a^{n-1}_j|$, for $n=1,\ldots,N$.  As in Fig. \ref{erra2}, the average values increase with $j$. Nonetheless, at the largest value, $j=64,000$, they are always smaller than $10^{-4}$. We can take this as an indication of the precision of the computed values $a_j$'s.

\begin{figure}
\centerline{\includegraphics[width=.6\textwidth, angle = -90]{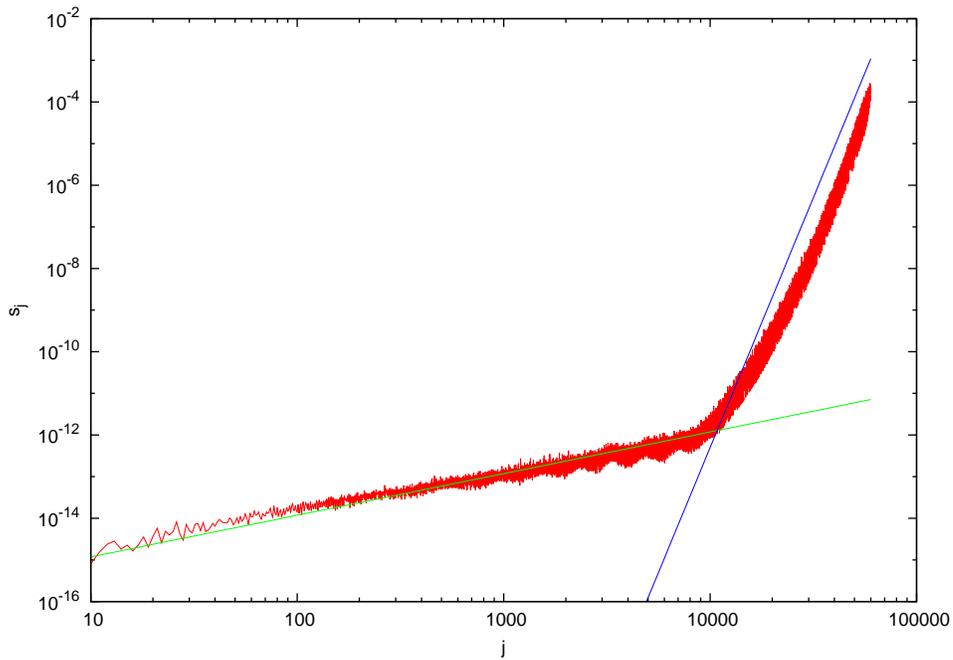}}
\caption{Standard deviation $s_j$ of the numerical error $\epsilon_j$ versus matrix index $j$  (red) and power--law law $s_j = A j$, with $A=1.18 \cdot  10^{-16}$ (green). A second power law, with exponent $B=12$ is also shown for comparison with the data for $j$ larger than 10,000 (blue).}
\label{expo0}
\end{figure}

A more refined analysis follows from considering the sequence $a_j(n)$ at convergence as a superposition of the true value $a_j$ and of a numerical error $\epsilon_j(n)$. Taking $\epsilon_j(n)$, for $n=1,\ldots,N$, as independent equally distributed random variables, we estimate their standard deviation $s_j$, that we take as a reliable measure of the error. Observe that we could have divided such standard deviation by the statistical factor $\sqrt{N-1}$, to yield the standard deviation of the sample average $\sum_{n=1}^N a_j(n) / N$. We prefer {\em not} to do this, because do not know whether the sample average converges to the true value $a_j$ and because we want to have a conservative estimate of the error of our technique.
Results are plotted in Figure \ref{expo0}, that shows an initial linear growth of the error, $s_j \sim A j$, which lasts until $j \simeq 10,000$, being later replaced by a more rapid increase, which appears to be faster than polynomial.

\begin{figure}
\centerline{\includegraphics[width=.6\textwidth, angle = -90]{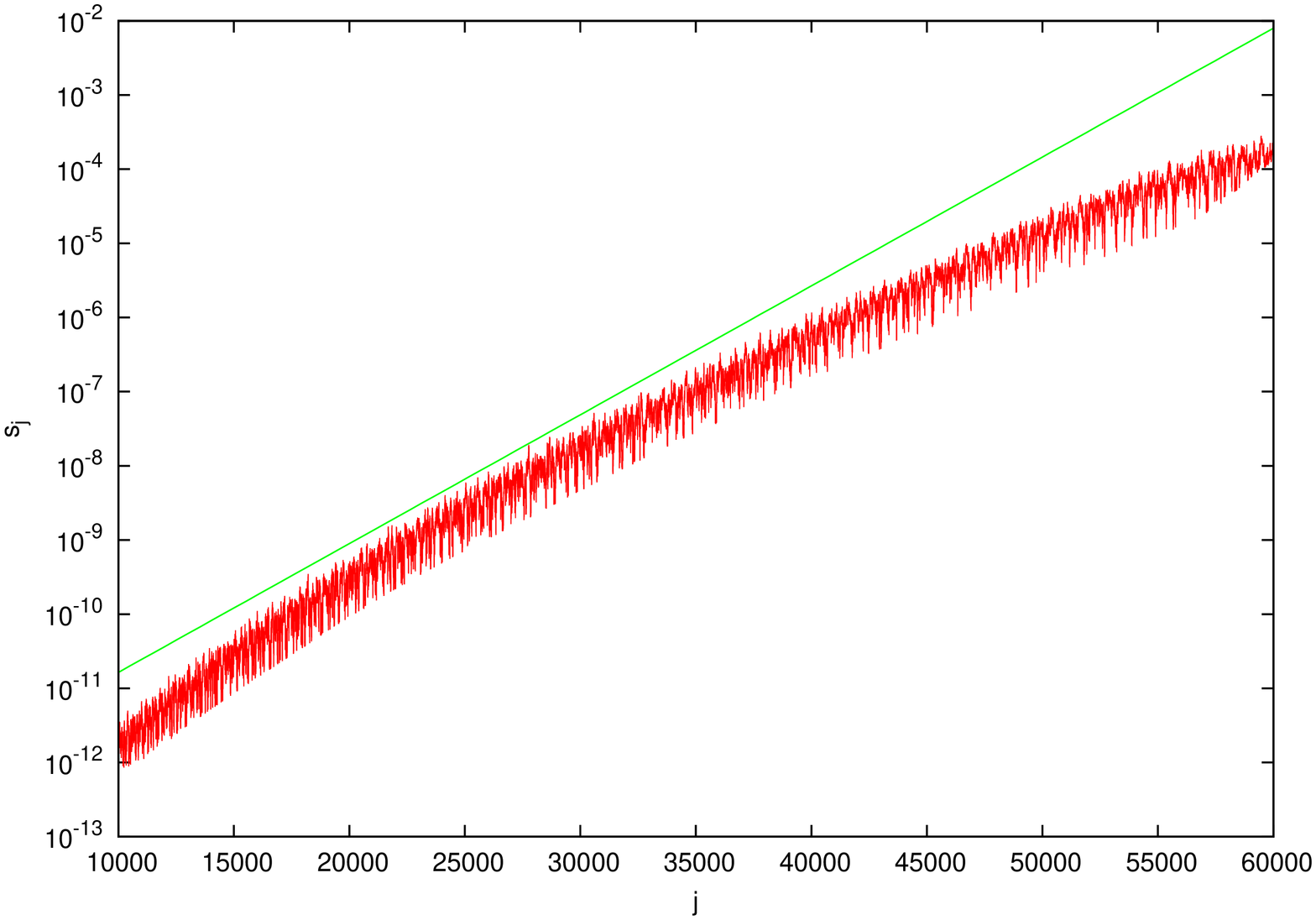}}
\caption{Standard deviation $s_j$ of the numerical error $\epsilon_j$ versus matrix index $j$  (red) and the exponential law $s_j = A e^{Bj}$, with $A=3 \cdot 10^{-13}$, $B=4 \cdot 10^{-4}$ (green).}
\label{expo}
\end{figure}

To investigate this second region we plot the data in semi--logarithmic scale, where a pure exponential appears as a straight line. This comparison shows that the error growth is less than exponential. The reasons beyond this behavior are not easily found, but we can nonetheless conclude that this technique allows us to compute reliably about 60,000 coefficients $a_j$, with the precision displayed in Figure \ref{expo0}.
An indirect, {\em a posteriori} proof of the precision of the computed Jacobi matrix is offered by the results of Section \ref{sec-numconf}.

\section{Regularity of Minkowski's measure}
\label{sec-mink}

\begin{figure}
\centerline{\includegraphics[width=.6\textwidth, angle = -90]{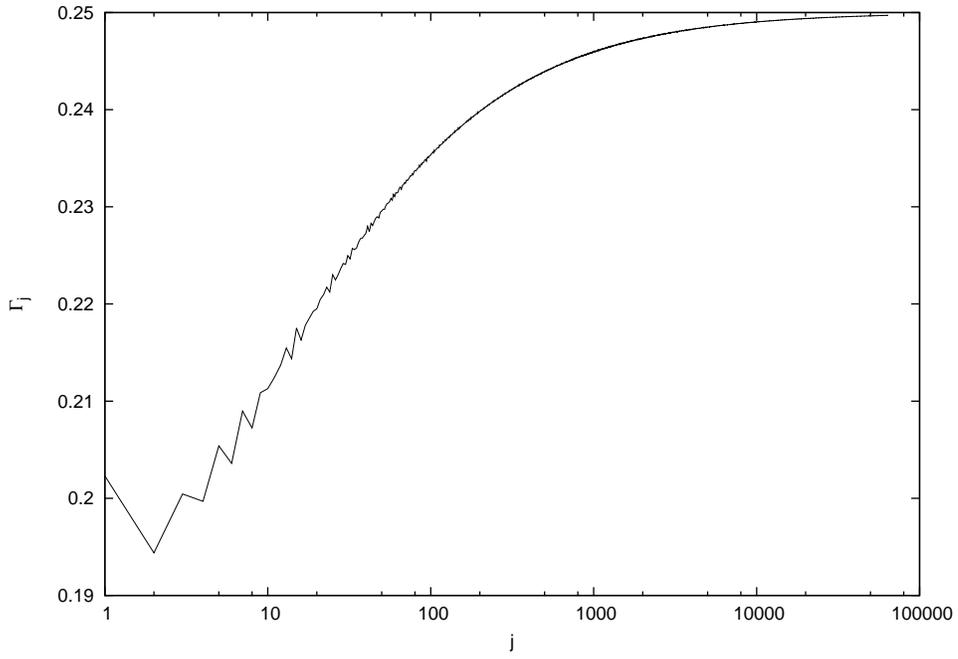}}
\caption{Geometric average $\Gamma_j$ in eq. (\ref{eq-reg01}) versus matrix index $j$.}
\label{fig-gamma1}
\end{figure}

We can now enter the heart of the matter.
We start by investigating the regularity of Minkowski's  measure $\mu$, defined as the existence of the limit
 $  \lim_{j \rightarrow \infty} \gamma_j^{-1/j} ={\mbox{Cap}([0,1])},$
where $\mbox{Cap}([0,1]) = \frac{1}{4}$ is the logarithmic capacity of the support of $\mu$. From eq. (\ref{nor2}) it easily follows that $1/\gamma_j = \prod_{l=1}^j a_l$; therefore, we need to compute the geometric average $\Gamma_j$ of the matrix entries $a_l$, for $l=1,\ldots,j$:
\beq
  \Gamma_j = \gamma_j^{-1/j} = \left[\prod_{l=1}^j a_l \right]^{1/j}.
  \label{eq-reg01}
  \end{equation}
We plot the results in Figure \ref{fig-gamma1}: we  observe that the conclusions of \cite{dresse} (convergence to a smaller value than $\frac{1}{4}$) were clearly due to the far--too--small range of $j$--values considered and to the deeper fact that the initial entries of the Jacobi matrix ``perceive'' a smaller support than $[0,1]$. We will explain this hitherto obscure remark in Section \ref{sec-zeroreg} and following. In the much wider range presented here convergence to the correct value appears.

\begin{figure}
\centerline{\includegraphics[width=.6\textwidth, angle = -90]{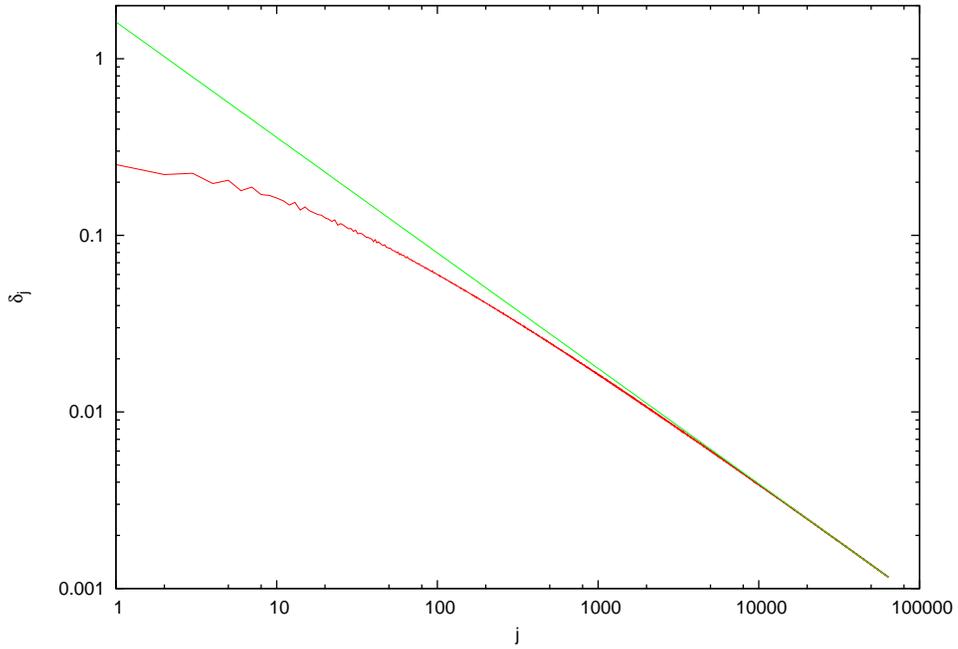}}
\caption{Difference $\delta_j$ in eq. (\ref{eq-reg02}) versus matrix index $j$ (red). Also plotted is the power-law decay $\delta_j = A j^{-B}$, with $A= 1.6186$, $B= 0.65424$ (green).}
\label{fig-gamma2}
\end{figure}

This conclusion can be put on even stronger grounds by studying the convergence rate. In fact, $\log(\Gamma_j)$ is the Ces\`aro average of the logarithm of the matrix entries $a_j$: the difference
\beq
  \delta_j = \log(\frac{1}{4}) - \frac{1}{j} \sum_{l=1}^j \log(a_l),
  \label{eq-reg02}
  \end{equation}
features a very clear asymptotic power-law decay of the kind $\delta_j = A j^{-B}$ (Figure \ref{fig-gamma2}). The pleasingly simple form of this law, its accuracy and the theoretical arguments of the following sections that explain the observed power--law type convergence from below, all confirm that the constant $\frac{1}{4}$ is the limit of $\Gamma_j$ when $j$ tends to infinity. We can therefore pretend that {\em the conjecture that measure is regular is proven numerically to a large degree of confidence}. We will add further evidence in favor of this conclusion in the following sections. Notice that having at our disposal matrix elements over various orders of magnitude of the index is essential.

Finally notice that the same data also show that $j \log 4 + \sum^j \log a_l$ does {\em not} converge to a finite value as $j$ tends to infinity. In Sect. \ref{sec-nevai} we will frame theoretically this result, relating it to the obvious lack of Szeg\"o asymptotics. We now turn to further verifications of the fact that the measure is regular.

\section{Zeros Minkowski's polynomials and their regularity}
\label{sec-zeroreg}

We now focus on regularity of the distribution of zeros of the polynomials, which, in our case, is equivalent to regularity of the measure. In this section and in the next we study whether the sequence $\nu_j$ in eq. (\ref{eq-nu3}) converges to $\nu_E$, the equilibrium measure on $[0,1]$, which is an absolutely continuous measure with density $\nu_E(x)$:
\beq
  \nu_E(x) = \frac{1}{\pi} \frac{1}{\sqrt{x(1-x)}}.
\nuq{eq-nu30}
We compute the zeros of $p_j(\mu;x)$ via diagonalization of the truncated Jacobi matrix $J(\mu)$. The algorithm (the implicit QL method) is notoriously stable and can be easily carried out to large orders. It yields the values $\zeta^j_l$, for $l=1,\ldots,j$, and $j$ in the range from 1 to 60,000.

To prove regularity, we check numerically a stronger property. Recall that the equilibrium measure $\nu_E$ is linked to the (properly shifted and normalized) Chebyshev polynomials, $p_j(\nu_E;x)$, whose zeros can be explicitly computed: Let
$\varphi^j_l= \frac{2l-1}{2j}$ for $l=1,\ldots,j$. Then,
\beq
  p_j(\nu_E;\theta^j_l) = 0 \; \mbox{   when  }  \theta^j_l = \frac{1}{2}[1-\cos(\varphi^j_l\pi)], \;\; l=1,\ldots,j.
\nuq{eq-nu31}
Numerical data support the validity of the following conjecture:
\begin{conjecture}
\label{conjzer}
As $j$ tends to infinity, the zeros $\zeta^j_l$ of the orthogonal polynomials of Minkowski's measure converge uniformly to the zeros of the Chebyshev polynomials of the same order. In addition, convergence is of power--law type: there exist positive constants $A$ and $B$ such that
\beq
  U_j = \max \{|\theta^j_l-\zeta^j_l|,  l=1,\ldots,j\} \leq A j^{-B} .
\nuq{eq-nu32}
\end{conjecture}
From the above, a simple consequence follows.
\begin{lemma}
\label{lem-reg}
Conjecture \ref{conjzer} implies regularity of the distribution of zeros.
\end{lemma}
{\em Proof}. Letting $f$ be a continuous function on $[0,1]$, for any $\epsilon >0 $ there exists $\delta>0$ such that $|\theta^j_l-\zeta^j_l| < \delta$ implies $|f(\theta^j_l)-f(\zeta^j_l)| < \epsilon$, where we use the trivial fact that continuity on a compact is necessarily uniform. Then, for any $\epsilon$ there exists $J$ such that
\beq
  |\frac{1}{j}\sum_{l=1}^j [f(\theta^j_l)- f(\zeta^j_l)]| \leq
     \frac{1}{j}\sum_{l=1}^j |f(\theta^j_l)- f(\zeta^j_l)| \leq \epsilon
\nuq{eq-nu33}
if $j > J$. Since the discrete measures generated by Chebyshev nodes converge weakly  to $\nu_E$, so does the sequence $\nu_j$. Observe that in this proof we do not need the estimate on the convergence rate. $\Box$

\begin{figure}
\centerline{\includegraphics[width=.6\textwidth, angle = -90]{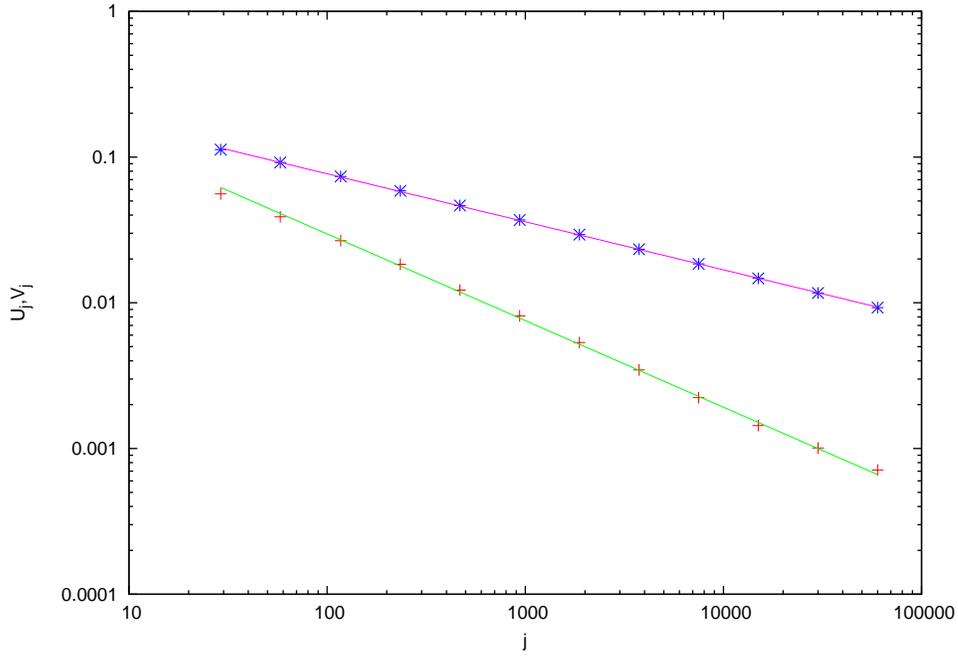}}
\caption{Upper bound $U_j$ from eq. (\ref{eq-nu32}) versus matrix index $j$ (crosses), fitted by the power-law decay $u_j = A j^{-B}$, with $A= 0.45$, $B= 0.59$. Upper bound $V_j$ from eq. (\ref{eq-zer1b}) (asterisks), fitted similarly, with $A=.35$, $B=0.32$}
\label{fig-tetzet}
\end{figure}

\begin{figure}
\centerline{\includegraphics[width=.6\textwidth, angle = -90]{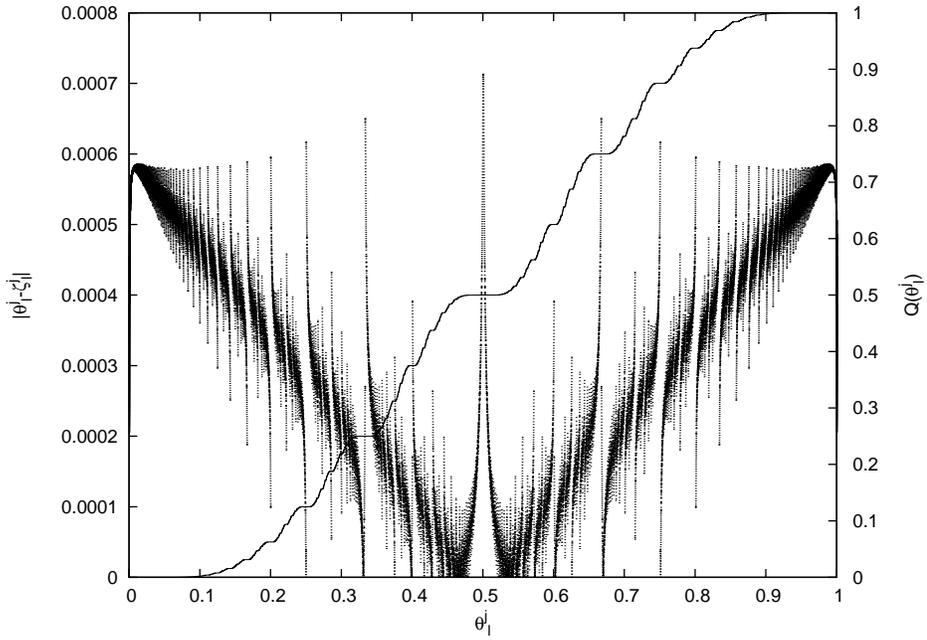}}
\caption{Absolute difference $|\theta^j_l-\zeta^j_l|$ versus $\theta^j_l$, for $j=60,000$ (left vertical scale). Each dot corresponds to a different value of $l$. Also plotted is the graph of the function $Q$ (continuous line, right vertical scale).}
\label{fig-tetzet2}
\end{figure}

Numerical support to Conjecture \ref{conjzer} is offered in Figure \ref{fig-tetzet}.
In Figure \ref{fig-tetzet2} we plot the absolute difference $|\theta^j_l-\zeta^j_l|$ versus $\theta^j_l$: this drawing reveals the locations where the two sequences differ the most. When compared to the graph of Minkowski's function $Q$, we observe that this pattern follows the structure of the Farey tree (defined in Remark \ref{rem-farey}): around rational values Minkowski's function has a non--analytic behavior, with large ``slippery plateaus'' ({\em i.e.} intervals of very low measure), which tend to ward off the zeros of $p_j(\mu;x)$. The phenomenon is stronger the higher the rational point is in the Farey tree. This fact is not surprising and will be further studied in the next sections.

To the contrary, the novel result of this investigation is that, this repulsion notwithstanding, Minkowski's measure is regular---said otherwise, zeros of the orthogonal polynomials find a way to ``creep in'' these regions and eventually populate them accordingly to the equilibrium measure. They do so at a power--law rate in the polynomial index, which explains the power--law behaviors found in Sections \ref{sec-mink} -- \ref{sec-discre}.
In Section \ref{sec-weights} we will encounter a dramatic consequence of these facts: the creation of {\em mathematical neutrinos}.

\section{Discrepancy Analysis}
\label{sec-discre}

Two measures $\lambda, \eta$ on $[0,1]$ differ as much as their integral over sub-intervals is different. This quantity is discrepancy:
\beq
 D(\lambda,\eta) = \sup \{ |\lambda(I)-\eta(I)|, \; I = (a,b) \subset [0,1] \}.
\nuq{eq-zer01}
Discrepancy has been intensely investigated \cite{andriev3,andriev4,er-turan}, since it is linked with root distribution and the logarithmic potential.

Consider the measures $\nu_j$ and $\nu_E$. If the sequence $D(\nu_j,\nu_E)$ converges to zero, then $\nu_j$ converges weakly to $\nu_E$ and this implies regularity of the zero distribution.
%
Zeros of the Chebyshev polynomials are best seen as projection on the real axis of equi-spaced points on the unit circle in the complex plane. When lifted to this set, the measure $\nu_E$ becomes the uniform Lebesgue measure. Therefore, it is convenient to also lift to the unit circle the zeros $\zeta^j_l$ of Minkowski's orthogonal polynomials: this leads to the definition of the normalized angles
\beq
\psi^j_l = \frac{1}{\pi} \arccos (1-2\zeta^j_l).
\nuq{eq-zer1}
The second set of data in Figure \ref{fig-tetzet} compares these values and the Chebyshev angles $\varphi^j_l$; the quantity $V_j$:
 \beq
 V_j = \max \{|\varphi^j_l-\psi^j_l|,  l=1,\ldots,j\},
 \nuq{eq-zer1b}
coherently with Conjecture \ref{conjzer}, it features a power--law decay.

Numerical computation of the discrepancy $D(\nu_j,\nu_E)$ is easily performed via a simple lemma that can be proven by straightforward computation:
\begin{lemma}
\label{lem-discrepa}
The discrepancy $D(\nu_j,\nu_E)$ between the discrete measure $\nu_j$
and the Chebyshev measure $\nu_E$ on $[0,1]$ is the maximum of the following three quantities:
\beq
D_1 = \max \{ |\psi^j_l - \frac{l - i}{j}|, \;  l=1,\ldots,j; \; i =0,1 \},
\nuq{eq-lem01}
\beq
D_2 = \max \{ |1-\psi^j_l - \frac{j-l+i}{j}|, \;  l=1,\ldots,j; \; i =0,1 \},
\nuq{eq-lem02}
\beq
D_3 = \max \{ |\psi^j_l-\psi^j_k - \frac{l-k+i}{j}|, \;  l,k=1,\ldots,j; \; i =-1,1 \}.
\nuq{eq-lem03}
\end{lemma}

\begin{figure}
\centerline{\includegraphics[width=.6\textwidth, angle = -90]{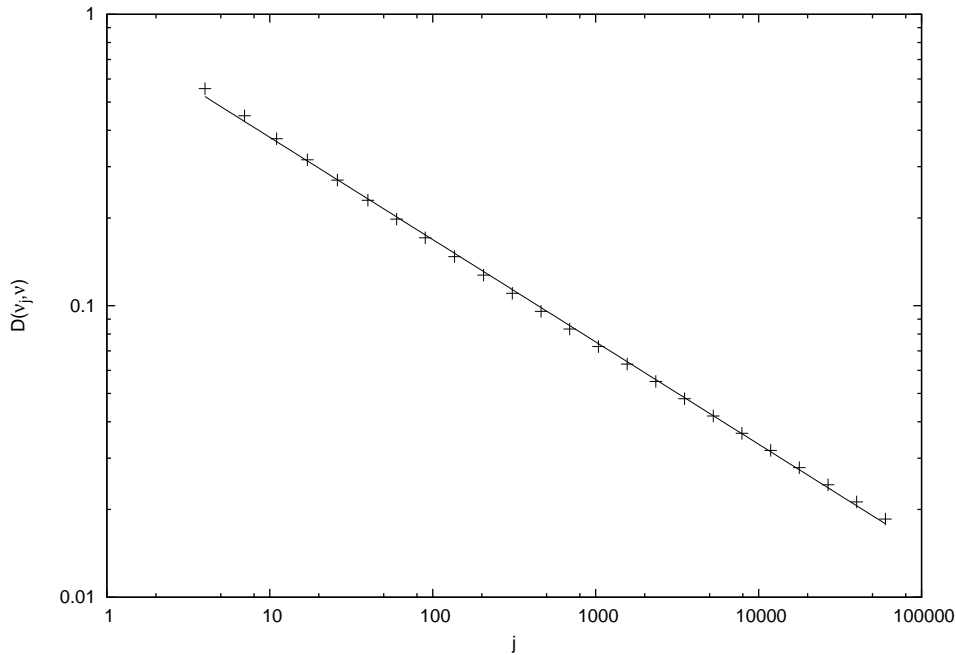}}
\caption{Discrepancy $D(\nu_j,\nu_E)$  versus $j$ and fit by the power-law decay $d_j = A j^{-B}$, with $A=.849$, $B= 0.351$.}
\label{fig-discrepa}
\end{figure}

Thanks to the lemma, numerical computation of $D(\nu_j,\nu_E)$ can be performed by a finite computation. This leads to the results plotted in Figure \ref{fig-discrepa}.
Again, not only we observe convergence towards zero of $D(\nu_j,\nu_E)$, but we can also extrapolate a power--law behavior with exponent close to $B=0.35$. Observe that by necessity, being $\nu_j$ discrete with atoms of equal weight, decay is bounded from below by $D(\nu_j,\nu_E) > 1/j$ for all $j$. We therefore formulate the following conjecture:

\begin{conjecture}
\label{conjdiscre}
As $j$ tends to infinity, the discrepancy $D(\nu_j,\nu_E)$ tends to zero. In addition, convergence is of power--law type: there exist positive constants $A$ and $B$ such that
$
   D(\nu_j,\nu_E)  \leq A j^{-B} .
$
\end{conjecture}

\section{Gaussian integration of Minkowski's  measure}
\label{sec-weights}

\begin{figure}
\centerline{\includegraphics[width=.6\textwidth, angle = -90]{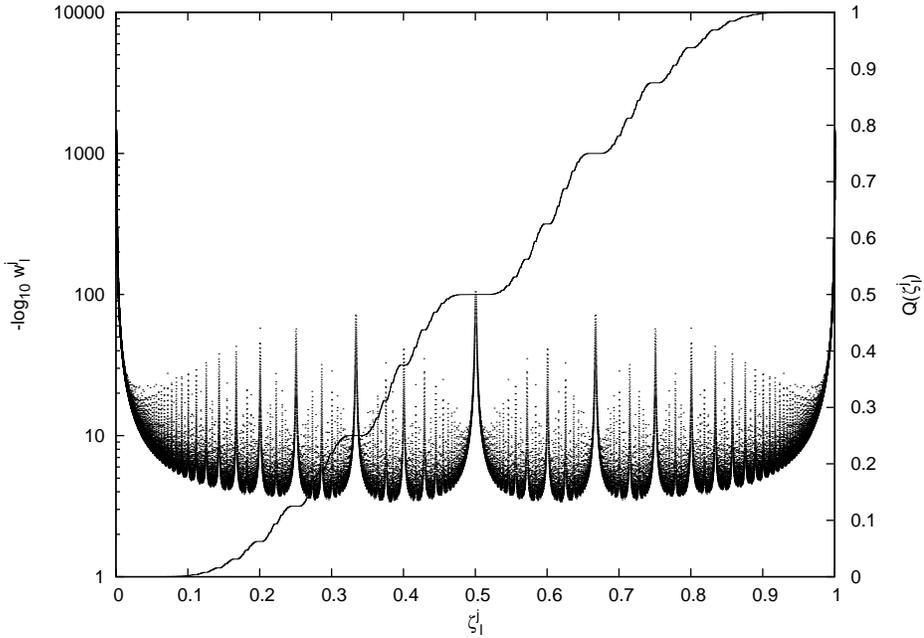}}
\caption{Base 10 logarithm of the inverse of the Gaussian weights $w^j_l$, plotted as dots at the abscissae $\zeta^j_l$, for $j=60,000$ (left vertical scale). Weights near $\zeta=0$ and $\zeta=1$ are smaller than $10^{-1000}$. Also plotted is the graph of the function $Q(\zeta)$ (continuous line, right vertical scale). The region near $x=\frac{1}{4}$ is magnified in Fig. \ref{fig-conve6bs} below.}
\label{fig-wei1b}
\end{figure}

We have seen that Minkowski's measure is {\em singular} (with respect to the Lebesgue measure) and at the same time {\em regular} (in the Ullman-Stahl-Totik sense). This hectic combination produces fascinating results, for what concerns Christoffel numbers---also called Gaussian weights. Recall eq. (\ref{eq-nu1}), which defines the discrete measure $\mu_j$, and let us focus on the weights $w^j_l$. It is convenient to start their analysis from a visual examination: in Figure \ref{fig-wei1b} we plot the {\em logarithm} of the inverse of the weights (in base 10) versus their location $\zeta^j_l$, for $j=60,000$. The figure is extremely instructive.

First of all, observe that the vertical scale at the left of the picture is also logarithmic, so that what appears graphically is minus {\em the logarithm of the logarithm} of the weight: the largest of these are of the order of $10^{-3}$, while the smallest are below $10^{-1000}$. There are about {\em one thousand orders of magnitude} of difference between the two!
We believe that it is then appropriate to call the latter {\em mathematical neutrinos}, for their elusive mass.
Their detection poses a serious challenge and requires a specific technique, by which it has been possible to compute Figure \ref{fig-wei1b}: this is explained in subsection \ref{sec-lyapu}.

An immediate application of these result is the computation of the Hausdorff dimension of the measure $\mu$, in Sect. \ref{sec-hausdim}.

Returning to Figure \ref{fig-wei1b}, observe that small weights appear when Gaussian nodes $\zeta^j_l$ approach the locations of rational points on the Farey tree (defined in Remark \ref{rem-farey}): this is the same to say, where Minkowski's $Q$ function (the continuous curve in the figure) appears to be ``almost flat''. Therefore,
Gaussian points and weights reflect the structure of Minkowski's measure and of the distribution of rationals on the real line \cite{paradis}, in a significant way. %
In subsection \ref{sec-slip} we develop an asymptotic theory which explicitly computes the ``cusps'' that appear in Figure \ref{fig-wei1b}.

Finally, in subsection \ref{sec-numconf} we deal with the fact that Farey points---the rationals---are dense in $[0,1]$: one must expect a dense set of these cusps. We show that their appearance  is hierarchically ordered, this fact being well described by the asymptotic theory that will be derived.

\subsection{Numerical computation of mathematical neutrinos}
\label{sec-lyapu}

Let again eq. (\ref{eq-nu1}) define the discrete measure $\mu_j$, with atoms at the Gaussian nodes $\zeta^j_l$. It is well known that they can be profitably computed as the eigenvalues of the truncated Jacobi matrix of rank $j$: $J^{(j)} u^j_l = \zeta^j_l u^j_l$. In the standard Golub-Welsch algorithm \cite{gene}, $w^j_l$ is obtained as the positive, first component of the normalized eigenvector $u^j_l$. The advantage of this algorithm lies in the fact that this component can be obtained without having to compute the full eigenvector. Unfortunately, this procedure is viable only as far as the amplitude of the weights does not fall below the level of numerical noise. When this happens---and Figure \ref{fig-wei1b} shows that this is the present case---the numerical detection of these weights poses challenges comparable to those of detecting neutrino's mass in physical experiments. In this subsection we describe a technique that can reveal these weights even when the ratio between the largest and the smallest is hundreds orders of magnitude smaller than machine precision.

Among many equivalent definition of the weights, we choose the one originally due to Shohat \cite{shoha,gautchr}.
Let $K_j(x,y)$ be the Christoffel-Darboux kernel:
\beq
   K_j(x,y) =   \sum_{l=0}^{j-1} p_l(\mu;x) p_l(\mu;y),
 \nuq{eq-nu2a}
The weights $w^j_l$ can be obtained from the diagonal values of this kernel at the roots of $p_j$ (that are computed as the eigenvalues of the truncated Jacobi matrix):
\beq
   (w^j_l)^{-1} = {K_j(\zeta^j_l,\zeta^j_l)} = { \sum_{l=0}^{j-1} p_l(\zeta^j_l)^2}.
 \nuq{eq-nu2}
It is customary to call the inverse of the diagonal kernels the Christoffel functions $\lambda_j(x)$:
\beq
   \lambda_j(x) =  \frac{1}{K_j(x,x)}.
 \nuq{eq-chr1}

The numerical calculation of these values is similar to the calculation of Lyapunov exponents in dynamical systems or in random matrix theory: according to eq. (\ref{eq-nu2a}), $K_j(x,x)$ is a summation of positive values. Low weights correspond to large values of the sum, which, in turn, arise from the fact that the norm of the two--components column vector $v_l = (p_l(\mu;x)$, $p_{l-1}(\mu;x))^t$ becomes very large and may overflow. The three--terms recurrence relation (\ref{nor2}) implies that $v_{l+1}$ is uniquely defined from the action of a {\em transfer matrix} acting upon $v_l$:
\beq
   v_{l+1} =
           \left(   \begin{array}{cc}  \frac{x-b_l}{a_{l+1}} &   - \frac{a_l}{a_{l+1}}    \cr
                               1 & 0    \cr
                               \end{array} \right) v_l.
 \nuq{eq-transf}
It is therefore not needed to store the full sequence $p_l(\mu;x)$ when calculating $K_j(x,x)$.
We then act as follows: starting from $l=1$, we compute the finite summation (\ref{eq-nu2}) term by term, until $K_l(x,x)$ exceeds a fixed threshold. At that moment, we let $V$ be the maximum of the absolute values of the two components of $v_l$, and we {\em renormalize} $v_l$ to $v_l/V$ and $K_l(x,x)$ to $K_l(x,x)/V^2$, while accumulating the logarithm of $V$ in a resummation variable $W(x)$, initially set to zero. At the end of the procedure, when $l=j-1$, the logarithm of the true value of $K_j(x,x)$ is obtained from the renormalized variables as $\log K_j(x,x) + 2 W(x)$ and the final formula follows:
\beq
   \log \lambda_j(x) =   - \log K_j(x,x) - 2 W(x).
 \nuq{eq-chr2}
Numerical values of Christoffel functions and numbers reported in this paper have been obtained using this procedure.

\subsection{Hausdorff dimension of the measure}
\label{sec-hausdim}

The first application of the precise numerical data of the previous subsection is the determination of {\em rigorous upper and lower bounds} to the Hausdorff dimension of Minkowski's measure $\mu$. Recall that this is defined as
$
\mbox{dim}_H(\mu) = \inf \{ \mbox{dim}_H (A), \mu(A)>0 \}
$, where the infimum is taken over all Borel measurable sets. It was proven by Kinney \cite{kinney} that one can write $\mbox{dim}_H(\mu)$ via an integral with respect to the measure $\mu$ itself:
\begin{equation}
\mbox{dim}_H(\mu) = \log(2) / [2 \int \log(1+x) d\mu(x) ].
\label{eq-dimh1}
\end{equation}
As a matter of facts, the function $f(x) = \log(1+x)$ in the integral is totally monotone on the positive real axis---{\em i.e.} $(-1)^n f^{(n)} (x) \geq 0$ on $[0,\infty)$ for all $n$. Under these conditions, the so-called first and second Gauss formulae of numerical integration give rigorous upper and lower bounds for the value of the integral in eq. (\ref{eq-dimh1}). The first formula is simply the integral of $f(x)$ with respect to the discrete measure $\mu_j$ in eq. (\ref{eq-nu1}): a finite summation, whose terms can be computed from Gaussian points and Christoffel weights. The second Gaussian formula is similar: it is based upon the truncated Jacobi matrix of the modified measure $d\rho(x) = x d\mu(x)$. This new Jacobi matrix $J(\rho)$ can be computed recursively starting from $J(\mu)$ \cite{gaut2}. Following these prescriptions we have obtained the data reported in Table \ref{tab1}. This permits to compute rigorously the Hausdorff dimension within 13 significative digits. These bounds contain the numerical value $\mbox{dim}_H(\mu) = 0.874716305108211142215152904219159757$ computed by Alkauskas \cite{alkatesi} via a convergent series involving the moments of $\mu$ in high precision arithmetics.

\begin{table}
\centering
\begin{tabular}{|l|l|l|l|}
  \hline
  $j$ & $\mbox{dim}_H(\mu)_+$ & $\mbox{dim}_H(\mu)_-$ & $\mbox{dim}_H(\mu)_+$ - $\mbox{dim}_H(\mu)_-$ \\
  \hline
   2 & 0.874761611261160  &0.874552879123086 & 0.000208732138074 \\
3 & 0.874716939422290 &0.874714034545017 & 0.000002904877273 \\
4 & 0.874716314143510 &0.874716274367535 &0.000000039775975 \\
5 & 0.874716305274063 & 0.874716304510136 & 0.000000000763927 \\
6 & 0.874716305110859 & 0.874716305099384 & 0.000000000011475 \\
7 & 0.874716305108267 & 0.874716305108003 & 0.000000000000264 \\
8 & 0.874716305108213 &0.874716305108207 & 0.000000000000006 \\
  \hline
\end{tabular}
\label{tab1}
\caption{Upper and lower bounds to the Hausdorff dimension of Minkowski's measure $\mu$ for increasing values of $j$, the number of points in Gaussian integration. The fourth column is the difference between the bounds.}
\end{table}

\subsection{Asymptotic analysis: Minkowski's slippery plateaus}
\label{sec-slip}
Let us now derive an asymptotic theory for the explanation of the Figure \ref{fig-wei1b}, which is based on the arithmetical nature of Minkowski's question mark measure, encoded by the M\'obius transformations (\ref{eq-mink1}).
In this subsection, we study the average value of Gaussian weights in certain subintervals $I \subset [0,1]$, when the polynomial degree $j$ tends to infinity. This average is defined by
\beq
 w^j_{I}  = \frac
 {\sum_{\zeta^j_l \in I} w^j_l}{\# \{l \mbox{ s.t. } \zeta^j_l \in I\} }.
\nuq{eq-neut01}
Since the vertical scale of Figure \ref{fig-wei1b} is logarithmic, we find appropriate to consider the logarithm of the average weight $w^j_{I}$. We start proving the following Lemma:
\begin{lemma}
\label{lem-3}
If Conjecture \ref{conj-reg} is verified,
\beq
 \log(w^j_{I})  = -\log (j) + \log(\mu(I) - \log(\nu(I)) + o(j),
\nuq{eq-neut04}
where $o(j)$ indicates as usual an infinitesimal sequence.
\end{lemma}
{\em Proof.}
Divide and multiply the denominator in eq. (\ref{eq-neut01}) by the order $j$, to obtain
\beq
j  w^j_{I}  = \frac
 {\sum_{\zeta^j_l \in I} w^j_l}{\# \{l \mbox{ s.t. } \zeta^j_l \in I\}/j }.
\nuq{eq-neut02a}
Because of weak convergence of the sequence $\mu_j$ to $\mu$, one has that the numerator at r.h.s. tends to the $\mu$--measure of $I$:
\beq
\lim_{j \rightarrow \infty} \sum_{\zeta^j_l \in I} w^j_l = \mu(I).
\nuq{eq-neut00}
Since the measure is regular (Conjecture \ref{conj-reg}), the denominator tends to $\nu(I)$. Taking logarithms yields $\lim \log(j) + \log  w^j_{I}  = \log(\mu(I)/\nu(I))$. $\Box$

We now specialize this result to particular intervals of the kind $[\frac{p}{q},\frac{p}{q}+y]$, in the neighborhood of the rational point $\frac{p}{q}$, where we observe numerically the cusps in the logarithm of the Christoffel numbers.
For sake of definiteness we consider the case $\frac{1}{q}$, which lies on the $(q-2)$-th level of the Farey tree (according to Remark \ref{rem-farey} we root the tree at $\frac{1}{2}$, level zero). Our technique can be applied also to the general case, with a minor additional effort. We need the following technical Lemma:

\begin{lemma}
\label{lem-farey}
The function $Q(x)$ is non-analytic around all Farey points $x=0$, $x=\frac{1}{q}$ with $q$ integer, in the following sense: for any integer $k \geq 0$ and $q>0$
\beq
 \mu([0,\frac{1}{k+1}]) = 2^{-k}; \;\;\;
 \mu([\frac{1}{q},\frac{k+1}{qk+q-1}]) = 2^{-k-q+1}.
\nuq{eq-broc90}

\end{lemma}
{\em Proof.}
Let $\Sigma$ be the set of finite words  in the letters $1$ and $2$. We denote by $|\sigma|$ the {\em length} of $\sigma \in \Sigma$: if  $|\sigma|=n$ then
$\sigma$ is the $n$-letters sequence $s_1,s_2,\ldots,s_n$ where $s_i$ is either $1$ or $2$. We associate to $\sigma$ the composite maps
\[
 M_\sigma = M_{s_1} \circ M_{s_2} \circ \cdots \circ M_{s_n}
 \]
 and
\[
 \Phi_\sigma = \Phi_{s_1} \circ \Phi_{s_2} \circ \cdots \circ \Phi_{s_n},
 \]
where the basic maps $M_i$ and $\Phi_i$ have been defined in eqs. (\ref{eq-mink1}) and (\ref{eq-mink11}).
Consider now the set $\Sigma_n$  of cardinality $2^n$ composed of all words of length $n$.
For any $n$, the graph ${\mathcal G}$ of $Q$ is contained in a union of images of the unit square:
\[
{\mathcal G} \subset \bigcup_{\sigma \; \mbox{\small  s.t.  } |\sigma|=n}  \Phi_\sigma([0,1]^2).
\]
In the above, each image $\Phi_\sigma([0,1]^2)$ is a rectangle of basis $M_\sigma([0,1])$ and height of length $2^{-|\sigma|}$. Since the graph ${\mathcal G}$ is continuous and since these rectangles are joined by the corners, it follows that the graph of $Q$ passes through these corners, so that
\beq
 \mu (M_\sigma([0,1])) = 2^{-|\sigma|}.
 \nuq{eq-broc01}
The basis intervals $M_\sigma([0,1])$ are the partition of $[0,1]$ produced by the Farey tree up to level $n-1$.

Consider now the word $\sigma=1^k$ composed of $k$ ones: we can use it to produce nested intervals shrinking at the point zero. In fact,
\beq
 M_{1^k}([0,1]) = [0,\frac{1}{k+1}], \;\;\; \mu (M_{1^k}([0,1])) = 2^{-k}.
 \nuq{eq-broc02}
This proves the first part of eq. (\ref{eq-broc90}).

We now consider the point $\frac{1}{q}$, with $q \geq 2$. Observe that $\frac{1}{2} = M_2(0)$ and that $\frac{1}{q} = M_{1}^{q-2} M_2 (0)$. This permits to map the intervals of the previous case, which were shrinking at the point zero, into a new sequence at the point $\frac{1}{q}$. Define $\sigma= 1^{q-2} 2 1^{k}$, to obtain
\beq
 M_{1^{q-2} 2 1^{k} }([0,1]) = [\frac{1}{q},\frac{k+1}{qk+q-1}], \;\;\; \mu (M_{1^{q-2} 2 1^{k}}([0,1])) = 2^{-k-q+1}.
 \nuq{eq-broc05}
The first part of the above equation can be proven by explicit calculation, while the second follows easily by computing the symbolic length of the word $\sigma$. This proves the second part of eq. (\ref{eq-broc90}).
$\Box$

Combining the previous results we have the following
\begin{lemma}
Let $q\geq 2$ and $k \geq 0$. Let $l_{q,k}=\frac{1}{q(qk+q-1)}$. Define the intervals \beq
I_{q,k}=[\frac{1}{q}+l_{q,k+1},\frac{1}{q}+l_{q,k}].
\nuq{eq-jqk}
The average amplitude of Gaussian points of order $j$ in $I_{q,k}$ behaves asymptotically as
 \beq
 \log (w^j_{I_{q,k}}) = -[(k+q) \log(2) + \log(j)   + \log(\nu(I_{q,k}))] + o(j).
\nuq{eq-neut04b}
\end{lemma}
{\em Proof.} Using the notations of Lemma \ref{lem-farey} the interval $I_{q,k}$ can be written as
\beq
 I_{q,k} = M_{1^{q-2} 2 1^{k} }([0,1]) \setminus M_{1^{q-2} 2 1^{k+1} }([0,1]),
 \nuq{eq-broc051}
and the first set at r.h.s. completely contains the second. To prove these facts, observe that the geometric length of the interval $M_{1^{q-2} 2 1^{k} }([0,1]) = [\frac{1}{q},\frac{k+1}{qk+q-1}]$ in eq. (\ref{eq-broc051}) is
\beq
  \frac{k+1}{qk+q-1}- \frac{1}{q} = \frac{1}{q(qk+q-1)} = l_{q,k}.
 \nuq{eq-broc06}
Using eq. (\ref{eq-broc05}) we then show that the measure of $I_{q,k}$ is $2^{-q-k}$. The thesis then follows using eq. (\ref{eq-neut04}). $\Box$

We can now state the principal result of this section:
\begin{proposition}
\label{prop-fundam}
Assume that Conjecture \ref{conj-reg} holds.
Let $0<y<\frac{1}{2}$, $q\geq 2$. Then, for sufficiently large $j$, the logarithmic amplitude $\log (w^j_{I_{q,k}})$, where $y \in I_{q,k}$, verifies the inequalities
\beq
- \log (w^j_{I_{q,k}}) \leq \log(2) [\frac{1}{q^2y} + \frac{1}{q} +q -2] + \log(j)
 - \frac{1}{2} H_-(q;y) + 2 \log (qy) - \log(1-q^2y);
 \nuq{eq-broc40}
 \beq
- \log (w^j_{I_{q,k}}) \geq \log(2) [\frac{1}{q^2y} + \frac{1}{q} +q - 3] + \log(j)
 - \frac{1}{2} H_+(q;y) + 2 \log (qy) - \log(1+q^2y),
 \nuq{eq-broc40b}
where $H_\pm(q;y)$ are continuous functions that tend to $\log(\frac{1}{q}-\frac{1}{q^2})$ when $y$ tends to zero.
\end{proposition}

{\em Proof.} Let $k$ be such that $y \in (l_{q,k+1},l_{q,k})$. This means that $\frac{1}{q}+y \in I_{q,k}$, the last set being defined in eq. (\ref{eq-jqk}) (see also eq. (\ref{eq-broc051})).
Since the measure $\nu_E$ is absolutely continuous, we use its density to write in the obvious way $\nu(I_{q,k}) =  |I_{q,k}| /  2  \sqrt{t_{q,k}-t_{q,k}^2}$, where $t_{q,k}$ is a point in  $I_{q,k}$ and where $|I_{q,k}|$ is the length of the interval $I_{q,k}$:
\[
  |I_{q,k}| = l_{q,k} - l_{q,k+1} = {q^2 l_{q,k} l_{q,k+1}}.
\]
Therefore, eq. (\ref{eq-neut04b}) becomes
\beq
 -\log (w^j_{I_{q,k}}) = (k+q-1) \log(2) + \log(j)  -\frac{1}{2} \log(t_{q,k}-t_{q,k}^2) +\log (q^2 l_{q,k} l_{q,k+1}) + o(j).
\nuq{eq-neut04c}
Since $k$ is such that $y \in (l_{q,k+1},l_{q,k})$, we have that
\beq
  \frac{1}{(k+1) q^2 + q(q-1)} \leq y \leq \frac{1}{k q^2 + q(q-1)},
 \nuq{eq-broc07}
and
\beq
 \frac{1}{q^2y} + \frac{1}{q} - 2 \leq k \leq
 \frac{1}{q^2y} + \frac{1}{q} -1.
 \nuq{eq-broc08}
Using this result, we can estimate the various quantities in eq. (\ref{eq-neut04c}).
\[
2 \log (qy) - \log(1+q^2 y) \leq
\log (q^2 l_{q,k} l_{q,k+1})  \leq
 2 \log (qy) - \log(1-q^2 y);
 \]
 \[
 t_{q,k} \leq \frac{1}{q} + l_{q,k} \leq \frac{1}{q} \frac{1+qy}{1-q^2 y};
 \]
 \[
 t_{q,k} \geq \frac{1}{q} + l_{q,k+1} \geq \frac{1}{q} \frac{1+qy}{1+q^2y}.
 \]
Let $h_\pm(q;y) = \frac{1+qy}{1 \mp q^2 y}$. These are bounded functions that tend to one as $y$ tends to zero. Consider the function $f(s)=s(1-s)$. This is a monotonic function, increasing for $s<\frac{1}{2}$ and decreasing for $s>\frac{1}{2}$. Therefore, the quantity $\log(t_{q,k}(1-t_{q,k}))$ can be estimated via the last two inequalities. We assume the first case, which corresponds to $q\geq 3$, the other consisting of reversed inequalities: in so doing, the definition
 \[
H_\pm(q;y)=\log[\frac{1}{q}h_\pm(q;y)(1-\frac{1}{q}h_\pm(q;y))]
\] completes the proof. $\Box$

\subsection{Hierarchical Christoffel functions of Minkowski's measure}
\label{sec-numconf}

Proposition \ref{prop-fundam} deals with the average value of Christoffel numbers over intervals in the neighborhood of the point $\frac{1}{q}+y$. As a matter of facts, only the value $y$ enters the formulae.  In addition, since $w^j_l$ is determined via eqs. (\ref{eq-nu2}), (\ref{eq-chr1}) we can think of $w^j_{J_k}$ as the typical logarithm of the Christoffel function at the point $\frac{1}{q}+y$. This permits to derive an asymptotic relation for this latter:

\begin{proposition}
\label{cor-chr}
Assume that Conjecture \ref{conj-reg} holds.
The typical logarithmic amplitude of the Christoffel function $\lambda_j(x)$ of order $j$ at the point $x=\frac{1}{q}+y$ is given by an asymptotic formula that comprises the sum of four contributions:
\beq
\log (\lambda_j(\frac{1}{q} + y))
\sim \Lambda_j(q;y) = \Lambda^1(q;y) + \Lambda^2(q) + \Lambda^3(q;y) + \Lambda^4(j),
 \nuq{eq-neut07}
in which the first is a geometrical factor that comes from the projection of the unit circle on the real axis:
\beq
 \Lambda^1(q;y) = \frac{1}{2} \log [\frac{1}{q} + y -(\frac{1}{q}+y)^2] + \log 2;
\nuq{eq-neut07a}
the second depends only on $q$, that is on the level in the Farey tree :
\beq
 \Lambda^2(q) = (-q + 2 - \frac{1}{q}) \log (2) + \log(\log(2));
\nuq{eq-neut07b}
the third contribution is determined by the distance from the rational point $\frac{1}{q}$ and it is singular in the limit of null distance:
\beq
 \Lambda^3(q;y) = - \frac{\log 2}{q^2 y} -2 \log (qy);
\nuq{eq-neut07c}
the fourth is minus the logarithm of the polynomial order:
\beq
 \Lambda^4(j) = - \log (j).
\nuq{eq-neut07c2}
\end{proposition}

{\em Proof.}
We use the results of Proposition \ref{prop-fundam}. Notice that the first term, $\Lambda_1(q;y)$,
can be written as $\log(|\sin \vartheta|)$, where $\vartheta  = \arccos [2(\frac{1}{q}+y-1)]$. In particular, it lies between the two bounds
\beq
  \frac{1}{2} H_-(q;y) + \log(2) \leq \Lambda_1(q;y) \leq \frac{1}{2} H_+(q;y) + \log 2.
   \nuq{eq-cor01}

The second contribution, $\Lambda_2(q)$, interpolates the purely $q$ dependent terms in eqs. (\ref{eq-broc40}),(\ref{eq-broc40b}). The choice of the interpolating constant $\log(\log(2))$ stems from a formal calculation of the asymptotics, not reproduced here. The estimate is justified by the inequalities
\[
 (-q -\frac{1}{q} + 1) \log(2) < \Lambda_2(q) < (-q -\frac{1}{q} + 2) \log(2)
 \]
Notice that in the bounds in the above equation differ from those of eqs. (\ref{eq-broc40}),(\ref{eq-broc40b}) by a single $\log(2)$ contribution, since this has been already included in the above eq. (\ref{eq-cor01}).

The third contribution, $\Lambda_3(q;y)$, collects the $y$ dependent quantities in eqs. (\ref{eq-broc40}),(\ref{eq-broc40b}) that diverge as $y$ tends to zero and therefore are at the origin of the cusps. Clearly, $\Lambda_3(q;y)$ interpolates the bounds because of the obvious inequalities $\log(1-q^2y) < 0 < \log(1+q^2 y)$. Furthermore, both of the latter functions are infinitesimal when $y$ tends to zero.

Finally, the last quantity, $\Lambda_4(j)$, directly follows by eqs. (\ref{eq-broc40}),(\ref{eq-broc40b}).
$\Box$


\begin{remark}
By using the symbolic approach of Lemma \ref{lem-farey}, similar results can also be derived for $\mu([\frac{p}{q},\frac{p}{q}+s])$, for any relatively prime $p$ and $q$, so that the local structure of the Christoffel function and of the diagonal Christoffel-Darboux kernel can be fully explained.
\end{remark}

Let us now illustrate the theoretical results with a series of figures, in which
we will use the symbol $K_j$, to underline the fact that the Christoffel functions are computed from eq. (\ref{eq-nu2a}) using the techniques of Section \ref{sec-lyapu}.

\begin{figure}
\centerline{\includegraphics[width=.6\textwidth, angle = -90]{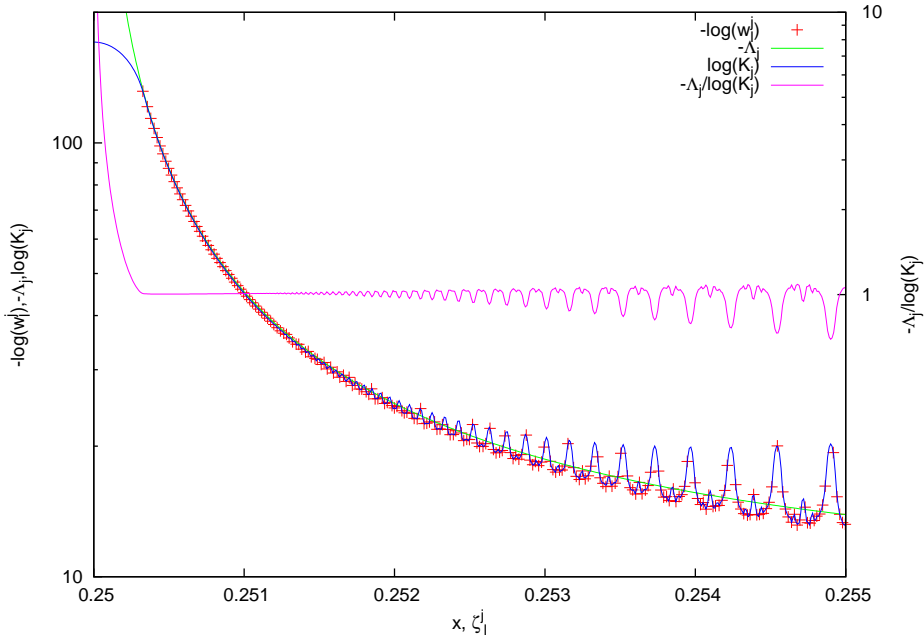}}
\caption{Logarithm of the Christoffel-Darboux kernel $K_j(x,x)$ in the right neighbourhood of $x=\frac{1}{4}$ (blue), with $j=60,000$, Gaussian pairs $(\zeta^j_l,-\log(w^j_l))$ (crosses) and asymptotic estimate $-\Lambda_j(\frac{1}{4};y)$ (with $y=x-\frac{1}{4}$) (green). The right vertical scale measures the ratio $-\Lambda_j(\frac{1}{4};y)/\log (K_j(x,x))$ (magenta).
See text for discussion.}
\label{fig-conve6bs}
\end{figure}

We start from Figure \ref{fig-conve6bs}, where we plot the logarithm of the diagonal kernel $K_j(x,x)$ for $x$ in the right neighborhood of the rational point $\frac{1}{4}$ (there is no particular reason behind this choice, other than it is a point on the second level of the Farey tree). The Gaussian pairs $(\zeta^j_l,-\log(w^j_l))$ lie on this graph, because of the fundamental property (\ref{eq-nu2}) and are plotted as crosses. Also, the asymptotic estimate $-\Lambda_j(\frac{1}{4};y)$ is drawn in the picture, where $y=x-\frac{1}{4}$ here and everywhere in the remainder of this section. The value of the polynomial order chosen is $j=60,000$.

Let us begin the analysis from the left part of the figure, {\em i.e.} $x$ values smaller than $0.2525$.
The distribution of Gaussian pairs reveals the detail of the cusp behavior already observed in Figure \ref{fig-wei1b}. The asymptotic formula $-\Lambda_j(\frac{1}{4};y)$ reproduces the cusp on which Gaussian points lie almost perfectly.
At the same time, $-\Lambda_j(\frac{1}{4};y)$ and $\log (K_j(x,x))$ are almost coincident as far as Gaussian pairs populate the graph and start parting at the location of the leftmost polynomial zero in the figure. In the innermost interval around $\frac{1}{4}$ the logarithm of the Christoffel kernel $K_j(x,x)$ stops following the cusp--like divergence: it is a continuous, bounded function and therefore it smoothly tends to a finite value when $x$ tends to $\frac{1}{4}$.

\begin{figure}
\centerline{\includegraphics[width=.6\textwidth, angle = -90]{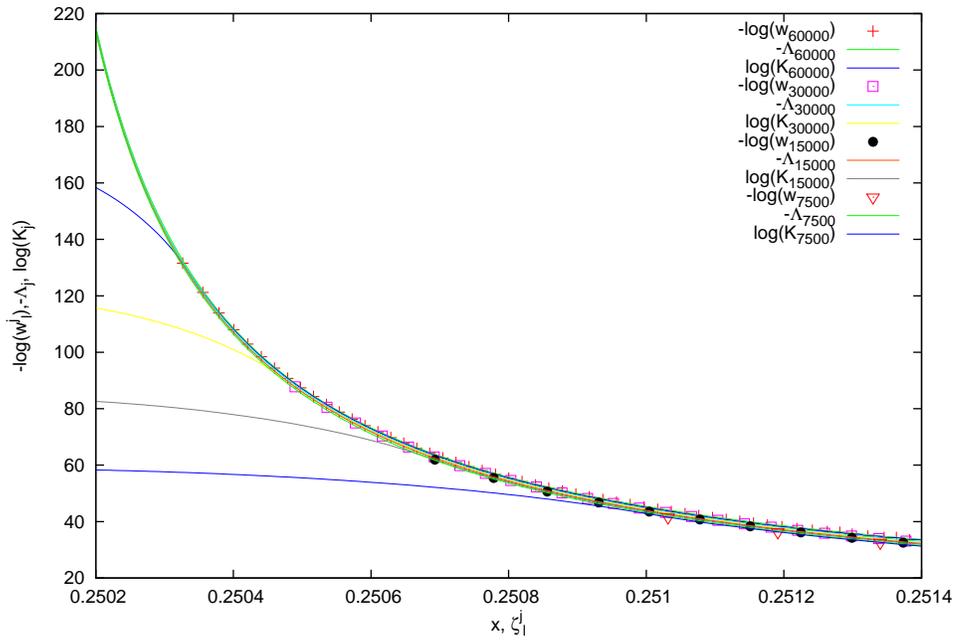}}
\caption{As in Figure \ref{fig-conve6bs}: logarithm of the Christoffel-Darboux kernel $K_j(x,x)$ in the right neighbourhood of $x=\frac{1}{4}$, Gaussian pairs $(\zeta^j_l,-\log(w^j_l))$ and asymptotic estimate $-\Lambda_j(\frac{1}{4};y)$ (with $y=x-\frac{1}{4}$). Geometrically increasing values of $j$ are plotted.}
\label{fig-conve1a6}
\end{figure}

These facts are all explained by the asymptotic theory. In fact, its derivation rests on the assumption that the measure is regular; since it describes the behavior of Gaussian weights, 
where there are no Gaussian points there is no reason to expect agreement with the kernel and one needs to wait for the asymptotic zero distribution to set in---see Sections \ref{sec-zeroreg}, \ref{sec-discre} and Figure
\ref{fig-conve1a6}.
The ratio between $-\Lambda_j(\frac{1}{4};y)$ and $\log (K_j(x,x))$ is also plotted in Figure \ref{fig-conve6bs}: it is divergent for $x=\frac{1}{4}$, but it rapidly approaches one when $x$ is equal to the leftmost polynomial zero, staying very close to this value in a certain interval. In the right part of the figure this ratio begins to oscillate while new peaks of the kernel appear: we will examine these momentarily, with the aid of a more detailed figure.

Let us now concentrate on the role of the polynomial order $j$. It enters the asymptotic relation (\ref{eq-neut07}) in a very simple way, via the term $\Lambda^4(j) = -\log(j)$. It neither determines the shape of the cusp, which is given by $\Lambda^3(q;y)$, nor it combines with the rational denominator $q$ or the location $x$. Also, when compared to other factors, the logarithmic dependence of $\Lambda^4(j)$ is rather mild. Yet, there is a hidden role of $j$ which does not appear explicitly 
and it is fully appreciated when considering the {\em asymptotic nature} of the expansion (\ref{eq-neut07}) in Proposition \ref{cor-chr}: {\em $\Lambda_j(\frac{1}{q};x-\frac{1}{q})$ approximates $\log \lambda_j(x)=-\log (K_j(x,x))$ better and better, in an interval which approaches $\frac{1}{q}$, as $j$ tends to infinity.}

In fact, let us first consider Figure \ref{fig-conve1a6}. It is analogous to Figure \ref{fig-conve6bs} but it plots quantities at geometrically increasing values of $j$, from $j=7,500$ to $j=60,000$, in a magnified range of values close to $\frac{1}{4}$. Zeros $\zeta^j_l$ are approximately equally spaced, with density given by the equilibrium measure $\nu_E$, in a region which excludes the immediate neighborhood of $\frac{1}{4}$. As before, the closest point $\zeta^j_l$ to the rational value $\frac{1}{4}$ marks the boundary of the region in which agreement between $\log(K(x,x))$ and $-\Lambda_j(\frac{1}{4};y)$ takes place. This range moves to the left when $j$ increases, according to Conjecture \ref{conjzer} in Section \ref{sec-zeroreg}.

\begin{figure}
\centerline{\includegraphics[width=.6\textwidth, angle = -90]{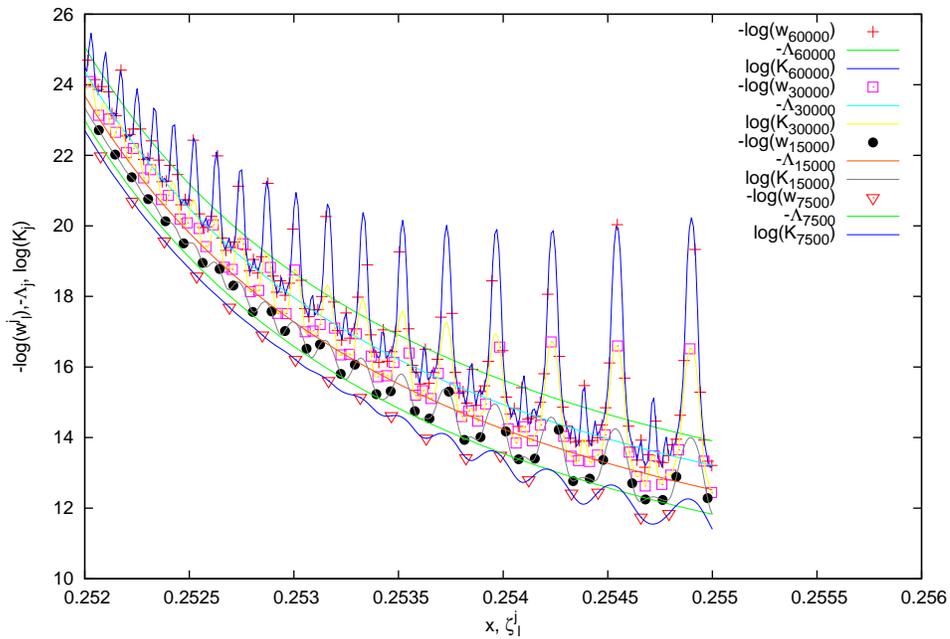}}
\caption{As in Figure \ref{fig-conve1a6}, for larger values of $x$. See text for analysis.}
\label{fig-conve3a}
\end{figure}

Secondly, eq. (\ref{eq-neut07}) 
holds {\em exactly} in the infinite $j$ limit for the {\em average} value of Christoffel numbers--Gaussian weights over the intervals $I_{q,k}$ (as follows from Proposition \ref{prop-fundam}), yet it is {\em asymptotic} for the {\em individual} weights, which are {\em pointwise} values of the Christoffel functions. This is seen in Figure \ref{fig-conve3a} that plots the magnification of the right part of Figure \ref{fig-conve6bs}. We observe the growth of new cusps superimposed to the leading  asymptotics $-\Lambda_j(\frac{1}{4};y)$, as the polynomial degree $j$ increases.
These cusps are associated with other rational points of higher order than $\frac{1}{q}$ in the Farey tree.
%
The same local asymptotic analysis leading to eq. (\ref{eq-neut07}) can be carried out at these new rational points, in a hierarchical construction that parallels that of Minkowski's question mark function by M\"obius IFS.

We can conclude that the behavior of Gaussian points--Christoffel numbers, described in this asymptotic analysis, reveals the hierarchical structure of Minkowski's measure in a very transparent way. This structure is also hierarchically encoded in the Jacobi matrix of the measure, in the sense that the finer structure the higher the index of the Jacobi matrix entries that reproduce it. This happens much in the same way as what observed in \cite{arx} for the Jacobi matrix of the equilibrium measure on fractal sets.

Finally, the fine details and the precision of the numerical results prove {\em a posteriori} that the procedures employed to compute Jacobi matrix, eigenvalues and Christoffel functions are stable and precise to the point of allowing such refined observations.

\section{The Nevai class of measures}
\label{sec-nevai}

\begin{figure}
\centerline{\includegraphics[width=.6\textwidth, angle = -90]{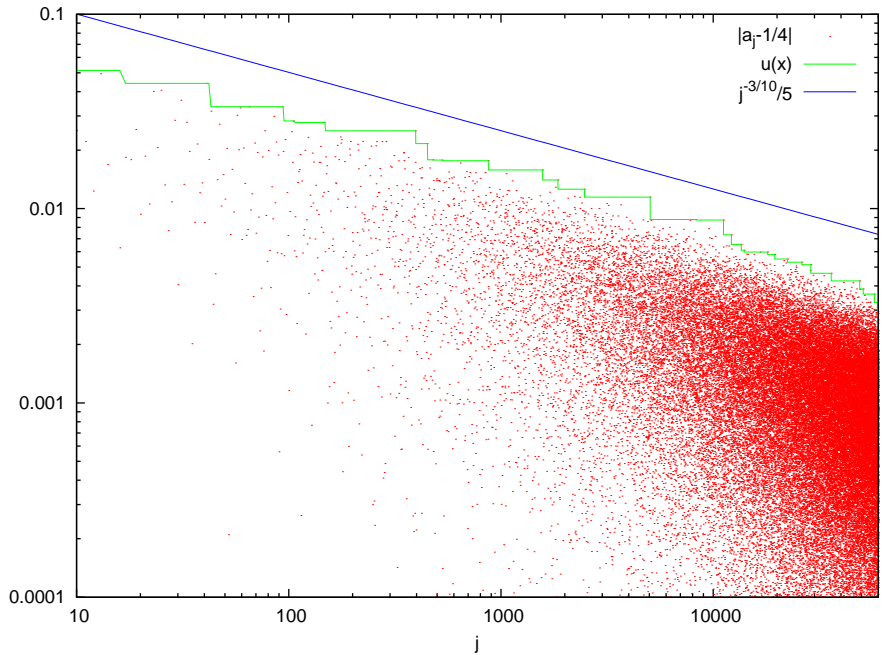}}
\caption{Convergence of the Jacobi matrix elements: $|a_j-\frac{1}{4}|$ versus $j$, the numerical upper bound $u(x)$ and the power law $j^{-3/10}/5$. Each red dot corresponds to a matrix element.}
\label{fig-aggei}
\end{figure}

The Nevai class of measures is defined by the fact that the matrix elements $a_j$ and $b_j$ tend to a limit. In our case, being all $b_j$'s equal to one half, this amounts to the limit $a_j \rightarrow \frac{1}{4}$.
In Figure \ref{fig-aggei} we plot the absolute difference $|a_j-\frac{1}{4}|$ in double logarithmic scale. Observe that the data, plotted with dots, range over orders of magnitude; yet, even the largest distances from the expected limit clearly tend to zero when $j$ grows. We also plot the numerical upper bound function $u(x)$, defined as $u(x)=\max\{|a_j-\frac{1}{4}|, x \leq j \leq 60,000\}$. This function decreases by definition, but it does so regularly. For comparison, we also plot a power-law decay: although with a lesser degree of confidence than all other numerical estimates in this paper, we may put forward the following
\begin{conjecture}
\label{conj-nevpower}
Convergence of the Jacobi matrix elements of Minkowski's measure is of power--law type: there exist two positive constants $A,B$ such that $|a_j-\frac{1}{4}| < A j^{-B}$.
\end{conjecture}
This convergence appears to be characterized by a considerably smaller exponent $B$ than the one observed in Figure \ref{fig-gamma2}. In fact, Figure \ref{fig-gamma2} is related to the average behavior of the logarithm of $a_j$ (regularity of the measure). To the contrary, the decay is here piloted by outliers that we need to master to prove the stronger requirements of the Nevai class.

In addition to the direct investigation of Jacobi matrix elements, a further analysis can be performed to test whether $\mu$ belongs to the Nevai class. Let us consider the ratio of the orthogonal polynomials $p_{j-1}(\mu;z)$ and $p_j(\mu;z)$. If $\mu$ belongs to $N(\frac{1}{4},\frac{1}{2})$ this ratio tends to the function $\phi_E(z) = 1/(z-\frac{1}{2} + \sqrt{z(z-1)})$, uniformly on compact subsets of the complement of $E=[0,1]$ (see the Introduction). At the same time, following {\em e.g.} \cite{walter2}, the above ratio can be written as the Stieltjes transform of a discrete measure $\sigma_j(\mu)$:
\begin{equation}
\frac{p_{j-1}(\mu;z)}{p_j(\mu;z)} = a_j \int \frac{d\sigma_j(\mu;s)}{z-s},
\label{eq-sig01}
\end{equation}
where, as in eq. (\ref{eq-sig02b}), we have put
\begin{equation}
\sigma_j(\mu) = \sum_{l=1}^n w^j_l  p_{j-1}^2(\mu;\zeta^j_l) \delta_{\zeta^j_l},
\label{eq-sig02}
\end{equation}
and, as in eq. (\ref{eq-nu2}), the Christoffel weights are
$(w^j_l)^{-1} = {\sum_{l=0}^{j-1} p_l(\mu;\zeta^j_l)^2}$. They obviously depend on $\mu$, even if this dependence is left implicit, not to overburden the notation.
We find convenient to introduce a short--hand notation for the weights in eq. (\ref{eq-sig02}):
\begin{equation}
S^j_l(\mu) =  w^j_l  p_{j-1}^2(\mu;\zeta^j_l).
\label{eq-sig02b2}
\end{equation}
Observe that $\sigma_j(\mu)$ is a sort of weighted version of the discrete measure $\mu_j$ associated with Gaussian integration, eq. (\ref{eq-nu1}). While the sequence of these latter {\em always} converges weakly (to $\mu$), convergence of $\sigma_j(\mu)$, when $j$ tends to infinity, is {\em not} always assured and we can use this fact for our purpose.

In fact, Nevai class $N(\frac{1}{4},\frac{1}{2})$, via convergence of the polynomial ratios to $\phi_E(z)$ and eq. (\ref{eq-sig01}) implies convergence of $\sigma_j(\mu)$ to the absolutely continuous measure $d\sigma_E(x) = \frac{1}{2 \pi} \sqrt{x(1-x)} dx$.
Viceversa, if $\sigma_j(\mu)$ converges weakly to a measure $\sigma$, then $\mu$ belongs to a Nevai class, since $\int x \; d\sigma_j(\mu;x) = b_{j-1}$ and $\int x^2 \; d\sigma_j(\mu;x) = b_{j-1}^2+a_{j-1}^2$. This proves the following criterion \cite{walter2}:
\begin{theorem}
The Minkowski's measure $\mu$ belongs to the Nevai class $N(\frac{1}{4},\frac{1}{2})$
if and only if the sequence of measures $\sigma_j(\mu)$ converges weakly to $\sigma_E$.
\end{theorem}

A first way to use this criterion comes from Conjecture \ref{conjzer} in Section \ref{sec-zeroreg}: we have noticed that zeros of the orthogonal polynomials $p_j(\mu;x)$ orderly converge to those of the Chebychev polynomials $p_j(\nu_E;x)$. It is immediate that $\nu_E$ belongs to the Nevai class $N(\frac{1}{4},\frac{1}{2})$ and therefore $\sigma_j(\nu_E)$ tends to $\sigma_E$. If we can control the distance between $\sigma_j(\mu)$ and $\sigma_j(\nu_E)$ we may expect to be able to prove convergence of $\sigma_j(\mu)$ to $\sigma_E$.

Indeed, let $\theta^j_l = \frac{1}{2}[1-\cos(\varphi^j_l\pi)]$ as in eq. (\ref{eq-nu31}) be the location of the roots of the Chebyshev polynomials $p_j(\nu_E;x)$. We can compute $\sigma_j(\nu_E)$ explicitly, since $p^2_{j-1}(\nu_E;\theta^j_l) = 2 \sin^2 (\frac{2l-1}{2n} \pi)$. This yields
\beq
 S^j_l(\nu_E) = 2 \sin^2 (\frac{2l-1}{2n} \pi) / n.
\nuq{eq-nucheb1}
Next, one might think of operating like in Lemma \ref{lem-reg}: let again $f$ be a continuous function,
\begin{eqnarray}
| \int f  d \sigma_j(\mu) - \int f d \sigma_j(\nu_E) | =
| \sum_{l=1}^j S^j_l(\mu) f(\zeta^j_l) - \sum_{l=1}^j S^j_l(\nu_E) f(\theta^j_l) |
\nonumber \\ \leq
\sum_{l=1}^j S^j_l(\nu_E) |f(\zeta^j_l) - f(\theta^j_l)| +
\sum_{l=1}^j  |f(\zeta^j_l)| |S^j_l(\nu_E)-S^j_l(\mu)| \nonumber  \\
 \leq
 \sup_l |f(\zeta^j_l) - f(\theta^j_l)|
 +
  \|f\|_\infty \sum_{l=1}^j   |S^j_l(\nu_E)-S^j_l(\mu)|.
\label{eq-nucheb3}
\end{eqnarray}
The first term at r.h.s. of the last inequality is infinitesimal, when $j$ tends to infinity,
according to Conjecture \ref{conjzer}. We can also verify numerically (see below, Fig. \ref{fig-sigma1}) that $s_j(\nu_E,\mu)=\max_l \{ |S^j_l(\nu_E)-S^j_l(\mu)| \}$ tends to zero as $j$ tends to infinity, and yet we find that the second summation at rhs, $\Sigma^0_j=\sum_{l=1}^j   |S^j_l(\nu_E)-S^j_l(\mu)|$ is not infinitesimal, but it seems to converge to a positive value: therefore, we cannot conclude from these estimates that $\sigma_j(\mu)$ converges to $\sigma_E$. We attribute this to the fact that the majorization made in eq. (\ref{eq-nucheb3}) is too crude, with respect to the fine properties of $\sigma_j(\mu)$.

We therefore need to compute numerically the Hutchinson distance $d(\sigma_j(\mu),\sigma_E)$. This can be done via an equivalent form of the definition (\ref{eq-hutdis}), which considers the distribution functions of the two measures: for any measure $\eta$ let $F(\eta;x) = \int \chi_{(0,x]}(y) \; d\eta(y)$. The Hutchinson distance becomes an integral of the absolute difference of the distribution functions:
\beq
  d(\sigma_j(\mu),\sigma_E) = \int_{0}^1 |F(\sigma_j(\mu);x) - F(\sigma_E;x)| \; dx.
\nuq{eq-hutdis2}
It is apparent that $F(\sigma_j(\mu);x)$ is a piece--wise constant function, which permits to split $[0,1]$ into intervals over which $F(\sigma_j(\mu);x)$ is constant and the difference $F(\sigma_j(\mu);x) - F(\sigma_E;x)$ has a fixed sign. The integral in eq. (\ref{eq-hutdis2}) so becomes a finite sum of integrals which can be explicitly computed in terms of elementary functions. The results obtained in this way are displayed in Figure \ref{fig-sigma1}. We observe a slow convergence towards zero of the distance $d(\sigma_j(\mu),\sigma_E)$ which is again dominated by a power--law:
\begin{conjecture}
\label{conj-nevsigma}
Convergence of the sequence of measures $\sigma_j(\mu)$ to $\sigma_E$ is of power--law type: there exist two positive constants $A,B$ such that the Hutchinson distance verifies
$ d(\sigma_j(\mu),\sigma_E) < A j^{-B}$
\end{conjecture}

\begin{figure}
\centerline{\includegraphics[width=.6\textwidth, angle = -90]{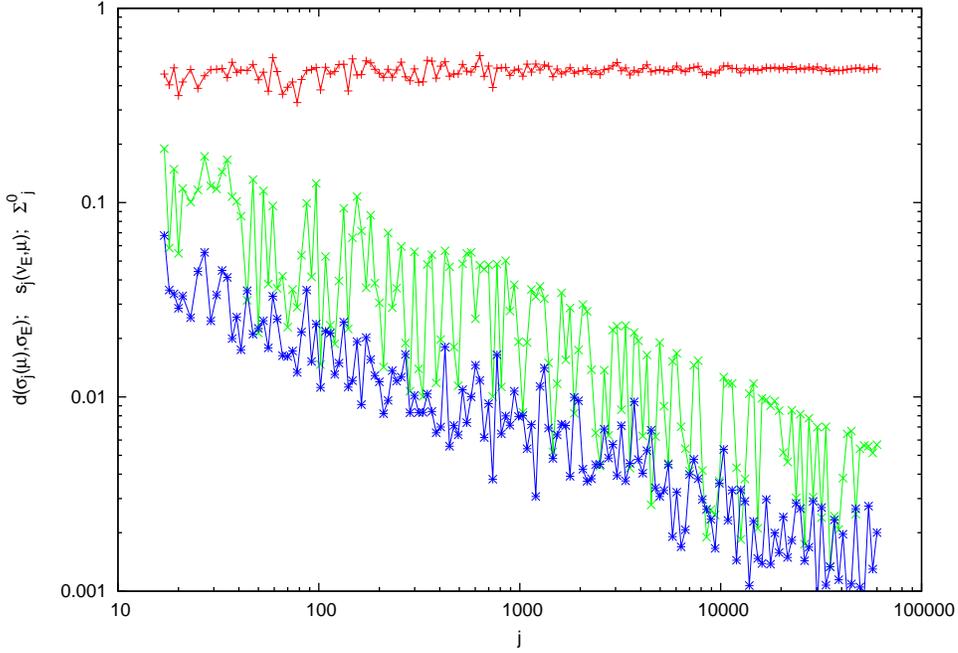}}
\caption{Quantities $\Sigma^0_j$ (red), $s_j(\nu_E,\mu)$ (green) and $d(\sigma_j(\mu),\sigma_E)$ (blue) versus matrix index $j$. Lines are merely to guide the eye and connect data at exponentially spaced values of $j$. See text for definition and discussion.}
\label{fig-sigma1}
\end{figure}

We conclude this work by showing that our numerical data are fully consistent with the theoretical fact that stronger convergence properties do not hold. For instance, power asymptotics, {\em i.e.} $\gamma_j/4^j \rightarrow \alpha$, with $0<\alpha<\infty$, can be easily seen to be equivalent to convergence of the series $\Sigma^3$, whose partial summation is defined, with other quantities to be discussed below, as follows:
\begin{eqnarray}
\Sigma^1_j = \sum_{l=1}^j | a_{l}-a_{l-1} |; \\
\Sigma^2_j = \sum_{l=1}^j |1- 16 \; a_{l}^2 |; \\
\Sigma^3_j = - \sum_{l=1}^j (\log (a_l) + \log(4)).
\end{eqnarray}

Convergence of the other series is regulated by the following ``collage'' theorem (reviewed {\em e.g.} in \cite{vanas,walter2}):
\begin{theorem}
If $\Sigma^2$ is finite, then the measure $\mu$ is absolutely continuous with respect to Lebesgue and its density belongs to Szeg\"o class.
Convergence of $\Sigma^2$ implies convergence of $\Sigma^1$. In turn, this latter implies that $\mu$ is absolutely continuous with respect to Lebesgue and its density is strictly positive and continuous on $(0,1)$.
\end{theorem}

\begin{figure}
\centerline{\includegraphics[width=.6\textwidth, angle = -90]{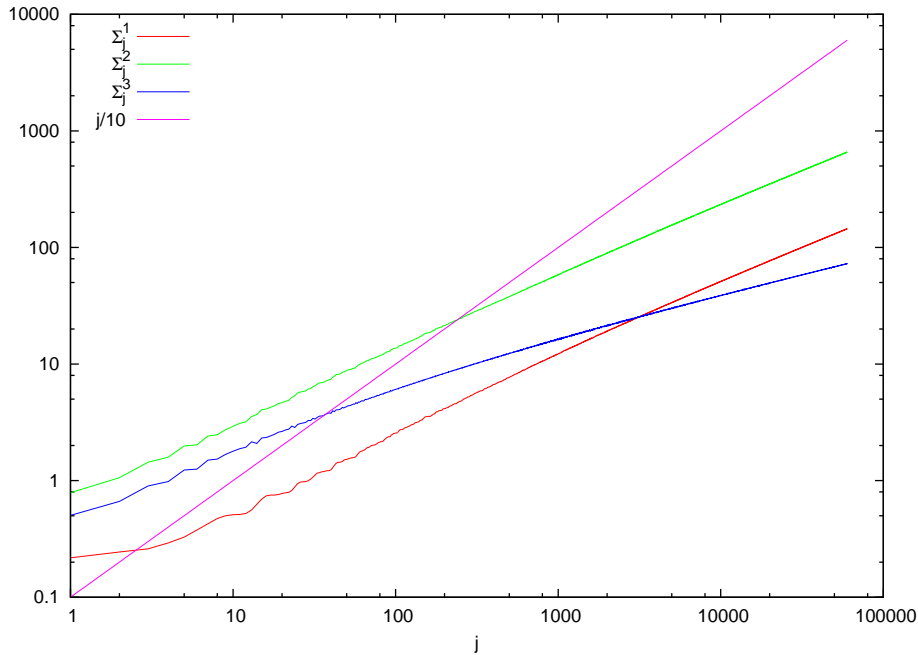}}
\caption{Quantities $\Sigma^m_j$ ($m=1,2,3$) versus matrix index $j$. See text for definition and comments.}
\label{fig-grow1}
\end{figure}

We expect divergence of the three series. This fact is clearly observed in Figure \ref{fig-grow1}. We also observe numerically that divergence of $\Sigma^3$ is of power--law type with exponent less than one. Since $\Sigma^3_j = j \delta_j$ (see eq. (\ref{eq-reg02}) and Fig. \ref{fig-gamma2}), this is consistent with regularity of the measure $\mu$.

\section{Conclusions}

As remarked by Totik \cite{totik2}, the Nevai class {\em seems to contain all sorts of measures}. In fact, it contains pure point measures \cite{alphwalt} as well as singular measures \cite{doron}. In addition, given any measure whose support is $[0,1]$, Totik has shown that there is a second measure, absolutely continuous with respect to the first, which {is} in Nevai class.
Therefore, it might seem that enlisting Minkowski's measure in this family is just another addition of minor interest. We think that this attitude is reductive, for two reasons. Firstly, since Minkowski's measure encodes the distribution of the rationals \cite{paradis} and their inner structure, proving that it does belong to the Nevai class might possibly reveal this structure from a different perspective. Secondly, no spectral characterization of Nevai class is known. Minkowski's measure falls short of verifying Rakhmanov sufficient condition \cite{rakh,totmany} {\em i.e.} almost everywhere positivity of the Radon Nikodyn derivative of $\mu$ with respect to Lebesgue---which is probably the widest sufficient condition known so far to this scope. Perhaps this fact is to be welcomed: Minkowski's measure could be the model of a possible widening of Rakhmanov condition. This investigation, as well as the simpler proof of regularity of Minkowski's measure, will be the object of future publications.


\end{document}